\documentclass[11pt,a4paper,reqno]{amsart}

\usepackage{amscd,amsmath}
\usepackage{amssymb,amsmath,latexsym,enumitem,tikz,dsfont,tikz-cd}
\usepackage[all]{xypic}       
\usepackage[varg]{pxfonts}
\usepackage{tabu,adjustbox,booktabs}
\usepackage{graphicx}
\usepackage{hyperref}
\usetikzlibrary{decorations.pathmorphing,shapes}

\usetikzlibrary{decorations.markings}
\tikzset{negated/.style={
        decoration={markings,
            mark= at position 0.5 with {
                \node[transform shape] (tempnode) {$\backslash$};
            }
        },
        postaction={decorate}
    }
}

\newcommand{\set}[1]{\{\,#1\,\}}

 \newcommand{\newdownarrow}{{{\rlap{$\ $}\hbox{$\downarrow$}}}}
 \newcommand{\newuparrow}{{{\rlap{$\ $}\hbox{$\uparrow$}}}}
 \newcommand{\twoheaddownarrow}{{\rlap{\rlap{$\ $}\raise .25ex\hbox{$\downarrow$}}\raise-.25ex\hbox{$\downarrow$}}}
 \newcommand{\twoheaduparrow}{{\rlap{\rlap{$\ $}\raise .25ex\hbox{$\uparrow$}}\raise-.25ex\hbox{$\uparrow$}}}

\newcommand{\tbigwedge}{\mathop{\textstyle \bigwedge }}
\newcommand{\tbigcap}{\mathop{\textstyle \bigcap }}
\newcommand{\tbigcup}{\mathop{\textstyle \bigcup }}
\newcommand{\tbigcoprod}{\mathop{\textstyle \coprod }}
\newcommand{\tbigprod}{\mathop{\textstyle \prod }}
\newcommand{\tbigvee}{\mathop{\textstyle \bigvee }}
 \newcommand{\cat}[1]{\ensuremath{\mathbf{#1}}} 
 
\makeatletter
\newcommand*{\@old@slash}{}\let\@old@slash\slash
\def\slash{\relax\ifmmode\delimiter"502F30E\mathopen{}\else\@old@slash\fi}
\makeatother

\def\Q{\mathbb{Q}}

\def\N{\mathbb{N}}
\def\Z{\mathbb{Z}}

\def\B{\mathcal{B}}

\def\SS{\mathsf{S}}

\def\FF{\mathcal{F}}
\def\EE{\mathcal{E}}
\def\K{\mathcal{K}}

\def\cl{\mathfrak{c}}
\def\op{\mathfrak{o}}
\def\bl{\mathfrak{b}}
\newtheorem{theorem}{Theorem}[section]
\newtheorem{proposition}[theorem]{Proposition}
\newtheorem{lemma}[theorem]{Lemma}
\newtheorem{corollary}[theorem]{Corollary}

\theoremstyle{definition}
\newtheorem{definition}[theorem]{Definition}

\theoremstyle{remark}
\newtheorem{remark}[theorem]{Remark}
\newtheorem{remarks}[theorem]{Remarks}

\author[I.~Arrieta]{Igor~Arrieta*}
\newcommand{\acr}{\newline\indent}

\address{\llap{*\,}Departamento de Matem\'aticas\acr
Universidad del Pa\'{\i}s Vasco UPV/EHU\acr
48080 Bilbao\acr
SPAIN}
\email{igor.arrieta@ehu.eus}
\begin{document}

\title{A dual Kuratowski--Mrówka theorem for locales}

\begin{abstract}The Kuratowski--Mrówka theorem states that a topological space or locale $X$ is compact precisely when every product projection $\pi_Y\colon X\times Y\to Y$ in the corresponding category is closed.

For topological spaces, the usual closure has a dual, referred to as  \emph{saturation}, which sends a subspace to the intersection of all open sets containing it. For locales, the corresponding operator is \emph{fitting}, which sends a sublocale to the intersection of all open sublocales containing it. 

We prove a dual form of the Kuratowski--Mrówka theorem for locales by
characterising compactness relative to fitting. For overt locales, we show constructively that this compactness property is equivalent to \emph{geometric semi-irreducibility}, a notion we introduce by requiring the locale classifying pairwise disjoint sequences of positive opens to be degenerate. Even classically, this result contrasts sharply with the situation in the category of topological spaces, where the analogous compactness for  saturation holds for every space.

Applied to discrete locales, our result yields a new characterisation of Freyd's notion of \emph{geometric finiteness} for objects in a topos. Together with the constructive Kuratowski--Mrówka theorem for ordinary localic compactness, this places {Kuratowski-finiteness} and geometric finiteness in a common framework as both arise from dual compactness properties.\end{abstract}
\maketitle

\section{Introduction}
In the category \cat{Top} of topological spaces and continuous maps, the usual closure has a dual construction, known as \emph{saturation}, obtained by interchanging the roles of closed and open subspaces. Whereas the closure of a subspace is the intersection of all the closed subspaces containing it, its saturation is the intersection of all open subspaces containing it. More generally, the theory of closure operators in categories finds applications in categorical topology
\cite{catfounds,DikranjanGiuli1987,Clementino1991,ClementinoGiuliTholen1996}.  In particular, this theory allows topological notions such as separation and compactness to be defined relative to a chosen closure operator. The usual closure and saturation provide two such choices in the category \cat{Top}.

For separation, the usual closure in \cat{Top} gives rise to the Hausdorff property, characterised by the closedness of the diagonal. With regard to compactness, the Kuratowski--Mr\'owka theorem \cite{Kuratowski1931,Mrowka1959} states that a topological space $X$ is compact if and only if the product projection $\pi_2\colon X\times Y\to Y$ in \cat{Top} is closed for every topological space $Y$. An easy equivalent formulation is that $\pi_2$ preserves the closures of subspaces.

In the category \cat{Loc} of locales, separation with respect to the usual localic closure is the so-called strong Hausdorff property (cf. \cite{separation}) whereas the Kuratowski--Mr\'owka theorem was extended by Pultr and Tozzi \cite{PultrTozzi1992} to the point-free setting and subsequently constructivised by Vermeulen \cite{propermaps}, who proved that a locale $L$ is compact if and only if the product projection $\pi_2\colon L\oplus 	M\to M$ in \cat{Loc} is closed for every locale $M$.

If one replaces closed sublocales with open sublocales, the point-free counterpart of the saturation operator is known as \emph{fitting} and sends each sublocale to the intersection of open sublocales containing it \cite{ClementinoPicadoPultr2018}. 

Separation with respect to fitting was first studied in \cite{APP21} under the name \emph{$\FF$-separatedness}, where it was shown that it has a strong structural parallel with the strong Hausdorff property: for instance, it is characterised by a condition of Dowker--Strauss type on the combinatorial structure of frame homomorphisms,  or it is characterised by certain weakened homomorphisms being frame homomorphisms. This duality was further explored by the author in \cite{Arrieta2024} in the context of other localic separation-type axioms.

In this setting, the question remains of identifying what compactness relative to the dual closure is. In \cat{Top}, the answer is easy: for all topological spaces $X$  and $Y$, the product projection $\pi_2\colon X\times Y\to Y$ sends saturated subspaces into saturated subspaces (see Proposition~\ref{saturation.compactness}). 

The situation is very different in the point-free context. We say a locale $L$ is \emph{$\FF$-compact} if every product projection $\pi_2\colon L\oplus M\to M$ sends fitted sublocales into fitted sublocales or, equivalently, if direct image along every such projection commutes with fitting. Our main result, Corollary~\ref{MAINRESULT}, characterises $\FF$-compactness for overt locales. An overt locale is $\FF$-compact if and only if the locale classifying the geometric theory of pairwise disjoint sequences of positive opens of $L$ is degenerate. We call this condition \emph{geometric semi-irreducibility}. Under classical logic, this says precisely that every pairwise disjoint family of nonzero opens is finite---a property introduced for topological spaces by A.~H.~Stone \cite{Stone1960} and subsequently extended  to locales by Dube \cite{irred}.

Interpreting our constructive results internally in a topos with a natural numbers object,  we obtain a new characterisation of Freyd’s notion of \emph{geometric finiteness} for objects in a topos. Following Freyd \cite{Freyd2006}, an object $X$ of a topos $\EE$ is {geometrically finite} if for any geometric morphism $f\colon \mathcal{B}\to\mathcal{E}$ with $\mathcal{B}$ a Boolean topos, the object $f^*(X)$ is Kuratowski-finite. We show that this is equivalent to $\FF$-compactness of the associated discrete locale $\mathcal{P}(X)$. Since $\mathcal{P}(X)$ is compact (in the usual sense) if and only if $X$ is Kuratowski-finite, Kuratowski-finiteness and geometric finiteness can be seen as duals to each other as they correspond to compactness relative to the usual closure and to fitting, respectively.

The paper is organised as follows. In Section~\ref{prelimFsep}, we provide preliminaries on locale theory, closure operators in a category and categorical topology. In Section~\ref{sect.dualcls} we introduce fitted localic maps and study their basic properties. In Section~\ref{sect4}, we use the locale of pairwise disjoint sequences to study semi-irreducibility in a constructive setting.  Section~\ref{sect.imp1}  and Section~\ref{sect.imp2} contain the proofs of the two directions of our characterisation of $\FF$-compactness. Section~\ref{section.discrete} specialises our result to discrete locales, where we establish the connection with Freyd's geometric finiteness. Finally, in Section~\ref{section.whcla} we explain the consequences of our main result under classical logic and  provide a simpler proof in that case.

The results in Sections~\ref{prelimFsep}--\ref{section.discrete} are constructively valid, except for those in Subsection~\ref{subsection.spaces} (which is about the category of topological spaces). Section~\ref{section.whcla} uses classical logic.

\section{Preliminaries}{\label{prelimFsep}}
\subsection{General remarks}
A set $X$ has \emph{decidable equality} if for all $x,y\in X$, either $x=y$ or $x\ne y$. We point out  that, crucially for this paper, the set $\N$ of natural numbers has decidable equality. 
Constructively, the usual notions of finiteness need not coincide. A set $X$ is \emph{Kuratowski-finite} if there is an $n\in \N$ and a surjective map $f\colon \{1,\dots,n\} \twoheadrightarrow X$ where $n=0$ is allowed; in other words, a set is Kuratowski-finite if its elements can be listed in a finite sequence, possibly with repetitions.  Kuratowski-finite sets have good closure properties; e.g., they are closed under Kuratowski-finite unions or under quotients. However, a subset of a Kuratowski-finite is in general not Kuratowski-finite. We recall  that a subset $A\subseteq X$ is \emph{complemented} if  for every $x\in X$ either $x\in A$ or $x\not\in A$. We shall need the following basic properties of  Kuratowski-finiteness.

\begin{lemma}\label{kf}The following properties hold:
\begin{enumerate}[label=\textup{(\arabic*)},leftmargin=2.0\parindent]
\item\label{kf1} A complemented subset of a Kuratowski-finite set is Kuratowski-finite;
\item\label{kf2} A Kuratowski-finite set is either empty or inhabited.
\end{enumerate}
\end{lemma}

 Given a set $X$, we denote by $\K(X)$ the \emph{Kuratowski-finite powerset of $X$}, namely the set of  the Kuratowski-finite subsets of $X$.

In Section~\ref{section.discrete}, we shall interpret constructive notions and results in the internal logic of toposes. All toposes considered in this paper are assumed to have a natural numbers object. For background on topos theory we refer to \cite{elephant}.

\subsection{Locale theory}\label{ssect.locale}
For constructive locale theory we refer the reader to \cite{Manuell2026} (see also  \cite[Chapter~C]{elephant}). Our notation and terminology regarding the categories of frames and locales will be that of \cite{PP12}; and for topics related to localic separation we refer to \cite{separation}. The Heyting operator in a frame $L$, right adjoint to the meet operator, will be denoted by $\to$; for each $a\in L$, $a^*=a\to 0$ is the pseudocomplement of $a$. Given a frame homomorphism $f\colon L\to M$, its right adjoint localic map is denoted by $f_*\colon M\to L$. Similarly, given a localic map $g\colon M\to L$, its left adjoint frame homomorphism is denoted by $g^*\colon L\to M$. 

If $\varphi$ is a proposition with no free variables in the ambient higher-order intuitionistic logic, we write $\llbracket \varphi\rrbracket$ for its truth value. We denote by $\Omega$ the set of truth values ordered by $\llbracket \varphi\rrbracket \leq \llbracket \psi\rrbracket$ if and only if $\varphi$ implies $\psi$; it contains the truth values $\top$ (true), the greatest element of $\Omega$, and $\bot$ (false),  the least element. With this ordering $\Omega$ is a frame, and in fact it is the initial frame.  For $S\subseteq\Omega$, joins can be computed as $\tbigvee S=\llbracket \top\in S\rrbracket$. A \emph{point} of a frame $L$ is a frame homomorphism $L\to\Omega$. We denote by $\mathsf{O}$ the \emph{trivial frame} with $0=1$, namely the final object in the category of frames.

\subsubsection*{Positivity} Constructively, asking that a set is nonempty is often not very useful; instead, one says that a set the set is \emph{inhabited} if it contains an element.  Similarly, given
a frame $L$, an element $a\in L$ is \emph{positive} if whenever $a\leq \tbigvee B$, then $B$ must me inhabited. Classically, this is equivalent to $a$ being nonzero. We write $a>0$ if $a\in L$ is positive, and moreover, we denote by
$$L^+:=\{ a\in L\mid a>0\}$$ the set of positive elements of $L$.

\subsubsection*{Overtness} A locale $L$ is \emph{overt} if the unique frame homomorphism $\Omega\to L$ has a left adjoint (classically, every locale is overt). There are several useful characterisations of overtness (cf. Proposition~\ref{overt.projec} below). For instance, $L$ is overt if and only if every element is a join of positive elements (in fact, if $L$ is overt and $A\subseteq L$, one has $\tbigvee A= \tbigvee\{ a\in A\mid a>0\}$). We shall need the following elementary property:

\begin{lemma}\label{property.overt} If $L$ is overt and $a\in L$, then $a$ is \emph{not} positive if and only if $a=0$. \end{lemma} \begin{proof}The `if' part is trivial so let us prove the `only if'.  Suppose that $a>0$ does not hold and, by overtness, express $a=\tbigvee B$ where $B\subseteq L^+$.  For every $b\in B$, one has $b\leq a$ and so $b>0$ is false (because $a>0$ is false). Hence $B=\varnothing$ and $a=0$. \end{proof}


\subsubsection*{Sublocales}\label{subsect.sublocales}

A \emph{sublocale} of a locale $L$ is a subset $S\subseteq L$ closed under arbitrary meets such that
\[\forall a\in L,\ \ \ \forall s\in S,\  \ \ a\to s\in S.\]
These are precisely the subsets of $L$ for which the embedding $j_S\colon S\hookrightarrow L$ is a morphism of locales. Sublocales of $L$ are in one-to-one correspondence with the regular subobjects (equivalently, extremal subobjects) of $L$ in \cat{Loc}. 

The system $\SS  (L)$ of all sublocales of $L$, partially ordered by inclusion, is a coframe \cite[Theorem~III.3.2.1]{PP12}, that is, its dual lattice is a frame.  Infima and suprema are given by
\[
\tbigwedge_{i\in I}S_i=\tbigcap_{i\in I}S_i, \quad \tbigvee_{i\in I}S_i=\set{\tbigwedge M\mid M\subseteq\tbigcup_{i\in I} S_i}.
\]
The least element is the sublocale $\mathsf{O}=\{1\}$ and the greatest element is the entire locale $L$. For any $a\in L$, the sublocales
\[
\mathfrak{c}_L(a)=\newuparrow  a=\set{b\in L\mid b\ge a}\ \text{ and }\ \mathfrak{o}_L(a)=\set{a\to b\mid b\in L}
\]
are the \emph{closed} and \emph{open} sublocales of $L$, respectively (that we shall denote simply by $\mathfrak{c}(a)$ and $\mathfrak{o}(a)$ when there is no danger of confusion). For each $a\in L$, $\mathfrak{c}(a)$ and $\mathfrak{o}(a)$ are
complements of each other in $\SS(L)$
and satisfy the expected identities
\begin{equation*}{\label{identities.basic}}
\tbigcap_{i\in I} \mathfrak{c}(a_i)=\cl(\tbigvee_{i\in I} a_i),\quad \cl(a)\vee\cl(b)=\cl(a\wedge b),
\end{equation*}
\[\tbigvee_{i\in I}\op(a_i)=\op(\tbigvee_{i\in I} a_i) \quad\mbox{ and }\quad \op(a)\cap \op(b)=\op(a\wedge b).
\]
Moreover, $\SS(L)$ is zero-dimensional; more precisely every sublocale is an intersection of sublocales of the form $\cl(a)\vee\op(b)$.

Given a sublocale $S$ of $L$, its \emph{closure}, denoted by $\overline{S}$, is the smallest closed sublocale containing it. In this context, the formula $\overline{S}=\cl(\tbigwedge S)$ holds. Moreover, $S$ is \emph{dense} if $\overline{S}=L$.

\subsubsection*{Fitted sublocales} A sublocale is \emph{fitted} if it is the intersection of all the open sublocales containing it. Hence, it constitutes a (non-conservative) point-free extension of the notion of saturated subspace. If $S$ is a sublocale  of $L$, the \emph{fitting} of $S$, denoted by $S^\circ$ is the intersection of the open sublocales containing $S$, that is, the smallest fitted sublocale containing $S$ (cf. \cite{ClementinoPicadoPultr2018}):
$$ S^\circ=\tbigcap_{S\subseteq\op(a)}\op(a).$$
We shall need the following easy property of fitted sublocales. In its proof, $\N$ can be replaced by any other index set with decidable equality; we prove it for $\N$ for our purposes:

\begin{lemma}\label{fitted.disjoint}
Let $L$ be a locale, $\{a_n\}_{n\in \N}\subseteq L$, and $\{S_n\}_{n\in \N}$ a family of sublocales satisfying the following properties:
\begin{enumerate}[label=\textup{(\arabic*)},leftmargin=2.0\parindent]
 \item\label{fitted.disjoint1} $S_n\subseteq \op(a_n)$ for each $n\in \N$, 
\item \label{fitted.disjoint2}$\op(a_n)\cap\op(a_m)=\mathsf{O}$ for all $n\ne m$,
\item \label{fitted.disjoint3}$S_n$ is fitted for each $n\in\N$. 
\end{enumerate}
Then, $\tbigvee_{n\in\N}S_n$ is fitted in $L$.
\end{lemma}

\begin{proof}
Let $S:=\tbigvee_{n\in\N} S_n$ and fix $n\in \N$. For any $a\in L$ with $S_n\subseteq \op(a)$, one has $S  \subseteq (\tbigvee_{m \ne n} \op(a_m))\vee \op(a)$ by \ref{fitted.disjoint1}, and consequently, $ S^\circ  \subseteq (\tbigvee_{m \ne n} \op(a_m))\vee \op(a)$. Hence, $ S^\circ \cap\op(a_n)  \subseteq  \op(a)$
by \ref{fitted.disjoint2}.  By fittedness of $S_n$, one concludes that $ S^\circ \cap\op(a_n)\subseteq S_n$. Finally, since joins of open sublocales distribute over any sublocale, one concludes 
$$ S^\circ = S^\circ  \cap \tbigvee_{n\in \N}\op(a_n)=\tbigvee_{n\in\N} \op(a_n)\cap  S^\circ \subseteq \tbigvee_{n\in\N} S_n=S.\qedhere$$
\end{proof}
\subsubsection*{Booleanization} Every locale $L$ has a largest dense sublocale denoted by $\mathcal{B}_L$; which is a complete Boolean algebra. Concretely it is given by
$$\mathcal{B}_L = \{ a\in L\mid a^{**}=a\} = \{ a^* \mid a\in L\}.$$

\subsubsection*{The Boolean algebra of smooth sublocales} Given a locale $L$, we shall need to consider the system $$\SS_b(L):=\bigl\{\tbigvee_{a\in A,\,b\in B} \cl(a)\cap \op(b) \mid A,B\subseteq L\bigr\}$$ ordered under sublocale inclusion.
\begin{lemma}[\cite{Arrieta2022}]\label{sbl}
For every locale $L$, the family $\SS_b(L)$ forms a complete Boolean algebra and there is a localic surjection $p_L\colon \SS_b(L)\twoheadrightarrow L$.
\end{lemma}
\subsubsection*{Image and preimage} Suppose that $g\colon L\to M$ is a localic map. 
Pullback along sublocales of $M$ yields a function $g_{-1}[-]\colon \SS(M)\to \SS(L)$, known as the \emph{inverse image of $g$}. It is well known that $g_{-1}[-]$ preserves arbitrary meets and (Kuratowski-)finite joins of sublocales---i.e., it is a coframe homomorphism, and satisfies the identities
$$g_{-1}[\op(a)]=\op(g^*(a))\quad\text{and}\quad g_{-1}[\cl(a)]=\cl(g^*(a)).$$ 
Furthermore, it has left adjoint colocalic map $$g[-]\colon \SS(L)\to \SS(M)$$
which takes every sublocale to its direct image (cf. \cite{PP12} for details). 
\begin{lemma}\label{lemma.factors}
Let $f\colon L\to M$ a frame homomorphism and $a\in L$. Then $f_*[M]\subseteq \op(a)$ if and only if $f(a)=1$.
\end{lemma}

\begin{proof}
One has $f_*[M]\subseteq \op(a)$ if and only if $M\subseteq f_{-1}[\op(a)]$ if and only if $M\subseteq \op(f(a))$ if and only if $a=1$.
\end{proof}
\subsubsection*{Open frame homomorphisms}
A frame homomorphism $f\colon L\to M$ is \emph{open} if  it has a left adjoint $f_!\colon M\to L$ satisfying the \emph{Frobenius identity}
\begin{equation}\label{fr}\tag{Fr}f_! ( f^*(b) \wedge a ) = b \wedge f_!(a)\end{equation}
for every $a\in  M$ and $b\in L$. A frame homomorphism $f\colon L\to M$ is open if and only if its right adjoint sends open sublocales into open sublocales---i.e., if $f_*[U]$ is open whenever $U$ is open. We shall need the following property:

\begin{lemma}\label{open.surj}
Let $f\colon L\to M$ be an open frame homomorphism. Then $f_*$ is surjective if and only if $f_!(1)=1$.
\end{lemma}

\begin{proof}
If $f_!(1)=1$, by \eqref{fr} one has $f_!(f^*(b))=f_! (f^*(b)\wedge 1)= b\wedge f_!(1)=b$ for any $b\in M$, hence $f_!$ is surjective. It follows that $f$ is injective and so $f_*$ is surjective by adjunction. The converse is trivial.
\end{proof}

In this context, we remark that a frame $L$ is overt if and only if the unique frame homomorphism $\Omega\to L$ is open. 

\subsubsection*{Coproducts of locales}  The coproduct of  $\{L_i\}_{i\in I}$ in \cat{Loc} will be denoted by $\tbigcoprod_{i\in I}^\cat{Loc} L_i$  in order to distinguish it notationally from the set-theoretic coproduct $\tbigcoprod_{i\in I}L_i$. Concretely, the frame of opens of $\tbigcoprod_{i\in I}^\cat{Loc} L_i$  is given by the set-theoretic cartesian product $\tbigprod_{i\in I}L_i$ ordered componentwise. In this setting, the left adjoint frame homomorphism to the coproduct injection  $$\iota_i \colon L_i \to \tbigcoprod_{i\in I}^\cat{Loc} L_i$$ 
  is simply the $i$-th projection.  If, moreover, the index set $I$ has decidable equality, the localic map $\iota_{i_0}$ can be described as t $\iota_{i_0}(a):=(a_i)_{i\in I}$ where $a_{i_0}=a$ and $a_i=1$ for $i\ne i_0$.

In particular, we shall denote the $I$-fold copower of $L$ in \cat{Loc} by $I{\cdot} L$---i.e.,
$$I{\cdot}L:=\tbigcoprod_{i\in I}^\cat{Loc} L.$$

\begin{lemma}\label{lemma.open.coproducts}
If $\{a_i\}_{i\in I}\subseteq L$, then the localic map $\tbigcoprod_{i\in I}^{\cat{Loc}}\op(a_i)\to L$ induced by the inclusions $\{\op(a_i)\hookrightarrow L\}_{i\in I}$  is open, and moreover it is surjective if and only if $\tbigvee_i a_i=1$.
\end{lemma}

\begin{proof}
The frame homomorphism $f\colon L\to  \tbigcoprod_{i\in I}^{\cat{Loc}}\op(a_i)$ given by $f(a)=( a_i \to a)_{i\in I}$ has a left adjoint given by $f_!( (x_i)_{i\in I} )= \tbigvee_{i\in I} a_i\wedge x_i$ satisfying \eqref{fr}. The second part follows from Lemma~\ref{open.surj}.
\end{proof}
\subsubsection*{Products of locales}  
Given locales $L$ and $M$ their product (equivalently, the coproduct of the corresponding frames) is denoted by $L\oplus M$ together with projections $\pi_1\colon L\oplus M\to L$ and $\pi_2\colon L\oplus M\to M$. The left adjoint frame homomorphisms, namely the coproduct injections, are denoted by $\iota_1\colon L\to L\oplus M$ and $\iota_2\colon M\to L\oplus M$. We denote
$$a\oplus b:= \iota_1(a)\wedge \iota_2(b),\quad\text{for }a\in L,\, b\in M.$$
In this setting, every element of $L\oplus M$ can be expressed as a join of elements of the form $a\oplus b$. For more details, we refer the reader to \cite[Section~2.4.1]{Manuell2026}.

A further characterisation of overtness can be given using product projections:

\begin{proposition}\label{overt.projec}
A locale $L$ is overt if and only if for every locale $M$ the projection $\pi_2\colon L\oplus M\to M$ is open.
\end{proposition}

We shall need the following description of the left adjoint of product projections (for a proof, see the proof of \cite[Corollary~3.25]{Manuell2026}):

\begin{lemma}\label{descript.lo.pi}
Let $L$ and $M$ be locales with $L$ overt and $a\in L^+$. Then $(\pi_2)_!(a\oplus b)=b$.
\end{lemma}

\subsection{Closure operators in a category}\label{sect.closop}

We start by recalling the definition of a closure operator in a category (in the sense of Dikranjan and Giuli \cite{DikranjanGiuli1987}, we also refer to \cite{Clementino1991,ClementinoGiuliTholen1996} and the references there for more information on closure operators). Given a factorisation system $(\mathcal{E},\mathcal{M})$, we denote by  $\mathsf{Sub}_{\mathcal{M}}(X)$ the partially ordered class of $\mathcal{M}$-subobjects of $X$.

\begin{definition}{\label{def.closure.op}}
Let $\mathcal{C}$ be a complete category equipped with a proper factorisation system $(\mathcal{E},\mathcal{M})$. A \emph{closure operator} $c$ with respect to $\mathcal{M}$ is  a family of maps
$$\{c_X\colon \mathsf{Sub}_{\mathcal{M}}(X)\longrightarrow \mathsf{Sub}_{\mathcal{M}}(X)\}_{X\in\mathsf{ob}(\mathcal{C})}$$
that satisfies the following properties:
\begin{enumerate}
\item it is monotone, that is, $m\leq n$ implies $c_X(m)\leq c_X(n)$ for all $m,n\in \mathsf{Sub}_{\mathcal{M}}(X)$;
\item it is inflationary, that is, $m\leq c_X(m)$ for all $m\in \mathsf{Sub}_{\mathcal{M}}(X)$;
\item for every $f\colon X\to Y$ in $\mathcal{C}$ and for all $m\in\mathsf{Sub}_{\mathcal{M}}(X)$, the identity $f(c_X(m))\leq c_Y(f(m))$ holds (where $f(m)$ refers to the $\mathcal{M}$-part of the $(\mathcal{E},\mathcal{M})$-factorisation of $f\circ m$).
\end{enumerate}
\end{definition}
A subobject $m\in \mathsf{Sub}_{\mathcal{M}}(X)$ is \emph{$c$-closed} $c_X(m)=m$ and   \emph{$c$-dense}  if $c_X(m)$ is an isomorphism. The closure operator $c$ is \emph{idempotent} if $c_X(m)$ is $c$-closed for every $m\in \mathsf{Sub}_{\mathcal{M}}(X)$ and every object $X$, and it is \emph{weakly hereditary} if the first map in the factorisation 
\[\begin{tikzcd}
	M & {c_X(M)} & X
	\arrow[from=1-1, to=1-2]
	\arrow["{c_X(m)}", from=1-2, to=1-3]
\end{tikzcd}\]
is $c$-dense for every $m\in \mathsf{Sub}_{\mathcal{M}}(X)$ and every object $X$. Finally, a morphism $f\colon X\to Y$ in $\mathcal{C}$ is \emph{$c$-preserving} if $f(c_X(m))\cong c_Y (f(m))$ for every $m\in \mathsf{Sub}_{\mathcal{M}}(X)$.

\subsubsection{$c$-separation}
An object $X$ in $\mathcal{C}$ is \emph{$c$-separated} if the diagonal morphism $\delta_X\colon X\to X\times X$ is $c$-closed. Separated objects enjoy useful categorical properties (cf. \cite[Proposition~4.2]{ClementinoGiuliTholen1996} and \cite[Propositions~10.1 and 10.7]{Clementino1991}). As is well known, for $\mathcal{C}=\cat{Top}$ the category of topological spaces with its usual surjection-embedding factorisation system and the usual Kuratowski-closure of subsets, the $c$-separated objects are precisely the Hausdorff spaces.

\subsubsection{$c$-compactness} An object $X$ in $\mathcal{C}$ is \emph{$c$-compact} if for every object $Y$ in $\mathcal{C}$, the product projection $\pi_2\colon X\times Y\to Y$ is $c$-preserving. The Kuratowski--Mr\'owka chararacterisation \cite{Kuratowski1931,Mrowka1959} states that a topological space $X$ is compact if and only if the product projection $\pi_2\colon X\times Y\to Y$ is closed for every topological space $Y$; hence, for $\mathcal{C}=\cat{Top}$ the category of topological spaces with its usual surjection-embedding factorisation system and the usual Kuratowski-closure of subsets, the $c$-compact objects are precisely the compact spaces.

In what follows we collect some of the main properties of $c$-compact objects:

\begin{proposition}[{\cite[Proposition~5.2]{ClementinoGiuliTholen1996}}]\label{properties.categorical.compactness}
 Let $\mathcal{C}$ be a complete category equipped with a proper factorisation system $(\mathcal{E},\mathcal{M})$ and a {closure operator} $c$ with respect to $\mathcal{M}$. Then the following properties hold:
\begin{enumerate}[label=\textup{(\arabic*)},leftmargin=2.0\parindent]
\item If $c$ is weakly hereditary and $m\in \mathsf{Sub}_{\mathcal{M}}(X)$ is $c$-closed and $X$ is $c$-compact, then $M$ is $c$-compact;
\item If $X$ is $c$-compact and $Y$ is $c$-separated, then every morphism $f\colon X\to Y$ is $c$-preserving;
\item\label{pr3} Let  $f\colon X\to Y$ is a morphism in $\EE$ such that every pullback of $f$ belongs to $\EE$. If  $X$ $c$-compact; then $Y$ is $c$-compact;
\item If $X$ and $Y$ are $c$-compact, then so is $X\times Y$.
\end{enumerate}
\end{proposition}

\begin{proof}
All the properties are stated in \cite[Proposition~5.2]{ClementinoGiuliTholen1996} except for \ref{pr3}, which follows from the proof in  \cite[Proposition~5.2]{ClementinoGiuliTholen1996} rather than from the statement.
\end{proof}

\section{The dual closure and fitted localic maps}\label{sect.dualcls}

We recall from Subsection~\ref{ssect.locale} that given a locale $L$ and a sublocale $S\subseteq L$,  the usual localic closure and the dual closure of $S$ are given by
$$\overline{S}=\tbigcap_{S\subseteq\cl(a)}\cl(a)\qquad\text{and}\qquad  S^\circ=\tbigcap_{S\subseteq\op(a)}\op(a),$$
---i.e., they define the smallest closed sublocale and the smallest fitted sublocale containing $S$ respectively.

As is well known, (Surjections, Regular monomorphisms) is a proper factorisation system in \cat{Loc}, and the assignment $S\mapsto \overline{S}$ defines a closure operator whose separated objects are the strongly Hausdorff ones (i.e., those with a closed diagonal, which have been extensively studied---see \cite{PP12,separation,catfounds}). Moreover, the compact objects with respect to this closure operator are precisely the compact locales---a result that goes back to Pultr and Tozzi \cite[Theorem~3.5]{PultrTozzi1992} and proved constructively by Vermeulen \cite[Proposition~3.7]{propermaps}.

On the other hand, it was shown in \cite[2.4]{ClementinoPicadoPultr2018}, among other properties of the dual closure, that the assignment $S\mapsto  S^\circ$ defines also a closure operator in \cat{Loc}, which will be referred to as \emph{fitting} operator. It is easily seen that this closure operator is idempotent and weakly hereditary (cf. Subsection~\ref{sect.closop}). 

\subsection{Fitted maps} A localic map $f\colon L\to M$ will be said to be \emph{fitted} if it is closure-preserving  with respect to the dual closure in \cat{Loc}---i.e., if $f[ S^\circ ]= f[S]^\circ$ for any sublocale $S\subseteq L$. The class of fitted localic maps is denoted by $\FF$. Next, we prove a characterisation of fittedness of localic maps similar to that of closed maps given in \cite[Lemma~1.5]{propermaps}---e.g., fitted maps are exactly those whose direct image preserves fitted sublocales.

Before proceeding, we introduce some notation. If $L$ is a locale, we denote by $\SS_f(L)$ the subset of $\SS(L)$ consisting of the fitted sublocales. Observe that $\SS_f(L)$ is a subcoframe of $\SS(L)$. Moreover, if $f\colon L\to M$ is a localic map, the restriction of the preimage yields a coframe homomorphism $f_{-1}[-]\colon \SS_f(M) \to \SS_f(L)$ whose left adjoint we denote by $f_{\dagger}[-]\colon \SS_f(L)\to \SS_f(M)$. Note that, more precisely,  $f_{\dagger}[S]=f[S]^\circ$ for any fitted $S\subseteq L$. With this preparation, we have:

\begin{proposition}
Let $f\colon L\to M$ be a localic map.  The following are equivalent:
\begin{enumerate}[label=\textup{(\arabic*)},leftmargin=2.0\parindent]
\item\label{fitcar1}  $f\in \FF$;
\item\label{fitcar2} $f$ sends fitted sublocales into fitted sublocales via direct image---i.e., $f[S]$ is fitted in $M$ whenever $S$ is fitted in $L$;
\item\label{fitcar3} $f_\dagger [ f_{-1}[U]\cap S]= U\cap f_{\dagger}[S]$ for all fitted $S\subseteq L$ and open  sublocale $U\subseteq M$.
\item\label{fitcar4} $(U\cap f[S])^\circ= U\cap {f[S]}^\circ$ for all fitted $S\subseteq L$ and open sublocale $U\subseteq M$.

\end{enumerate}
\end{proposition}

\begin{proof}
First, the equivalence between \ref{fitcar1} and \ref{fitcar2} follows from the idempotency of the closure operator (cf. \cite[3.1]{ClementinoGiuliTholen1996}). 
Note also that the equivalence between \ref{fitcar3} and \ref{fitcar4} follows at once by using the explicit description of $f_\dagger[-]$, \cite[Proposition~1.4]{propermaps} and the fact that open sublocales are complemented. Moreover, that \ref{fitcar2} implies \ref{fitcar4} is trivial. For proving that \ref{fitcar4} implies \ref{fitcar2}, let $S$ be fitted in $L$. We need to show that $f[S]^\circ\subseteq f[S]$. By  zero-dimensionality of $\SS(L)$ it suffices to show that if $F$ is closed and $U$ is open in $M$, $f[S]\subseteq U\vee F$ implies $f[S]^\circº\subseteq  U\vee F$.  This is clearly equivalent to showing that for any open sublocales, $U,V\subseteq M$, $U\cap f[S]\subseteq V$ implies $U\cap f[S]^\circ\subseteq V$. Indeed, since $U\cap f[S]\subseteq V$ implies $(U\cap f[S])^\circ \subseteq V$,  the conclusion follows from \ref{fitcar4}.
\end{proof}

\subsection{$\FF$-separation} Separated objects in \cat{Loc} with respect to the fitting operator were first studied in \cite{APP21} under the name \emph{$\FF$-separated locales}. This separation property is weaker than fitness and exhibits a remarkable duality with the strong Hausdorff property. This duality was subsequently studied and developed further in \cite{Arrieta2024}.

\subsection{$\FF$-compactness} 
We now turn into the main topic of the paper. A locale $L$ will be said to be \emph{$\FF$-compact} if it is compact with respect to the fitting operator---i.e., if the product projection 
$$\pi_2\colon L\oplus M\to M$$
is fitted for every locale $M$ (that is, if $\pi_2[S]$ is fitted in $M$ whenever $S$ is fitted in $L\oplus M$). The following properties follow at once from the  categorical setting in Proposition~\ref{properties.categorical.compactness}:

\begin{proposition}\label{properties.categorical.compactness.fitted}
  Let $L$ and $M$ be locales. Then the following properties hold:
\begin{enumerate}[label=\textup{(\arabic*)},leftmargin=2.0\parindent]
\item\label{properties.categorical.compactness.fitted1} If $S$ is fitted in $L$ and $L$ is $\FF$-compact, then $S$ is $\FF$-compact;
\item\label{properties.categorical.compactness.fitted2} If $L$ is $\FF$-compact and $M$ is $\FF$-separated, then every localic map $f\colon L\to M$ is fitted;
\item\label{properties.categorical.compactness.fitted3} Let  $f\colon L\to M$ a localic surjection such that every pullback of $f$ is surjective. If  $L$ is $\FF$-compact, then $M$ is $\FF$-compact;
\item\label{properties.categorical.compactness.fitted4} If $L$ and $M$ are $\FF$-compact, then so is $L\oplus M$.
\end{enumerate}
\end{proposition}

\subsection{The saturation operator in the category of topological spaces}\label{subsection.spaces}

The study of $\FF$-compactness in $\cat{Loc}$ naturally suggests an analogous question in the category $\cat{Top}$ of topological spaces and continuous maps. Recall that the counterpart to the fitting operator in  \cat{Top} is the \emph{saturation operator} sending a subspace to the intersection of all open subspaces containing it. A closure-preserving map with respect to the saturation operator is simply a  continuous map   which maps saturated subspaces into saturated subspaces.  As is well known, a topological space is separated with respect to the saturation operator if and only if it is $T_1$.

In view of the discussion above concerning $\FF$-compactness one might wonder what the compact spaces are with respect to the saturation operator. We answer this topological question below, showing that every space satisfies the condition.
First, since (assuming classical logic) a power set is completely distributive, we have:

\begin{lemma}
Saturated subspaces are closed under arbitrary unions in \cat{Top}.
\end{lemma}

\begin{proposition}\label{saturation.compactness}
Every topological space is compact with respect to the saturation operator.
\end{proposition}

\begin{proof}
Let $X$ and  $Y$ be topological spaces. We have to show that $\pi_2\colon X\times Y\to Y$ sends saturated subspaces into saturated subspaces. Let $A=\tbigcap W_i$ be saturated in $X\times Y$ (where each $W_i$ is open). Since $A=A\cap \tbigcup_{x\in X}\{x\}\times Y=\tbigcup_{x\in X}A\cap \{x\}\times Y$, by the previous lemma and the fact that $\pi_2[-]$ preserves unions, it suffices to show that for each $x\in X$, the subspace $\pi_2[ A\cap \{x\}\times Y ]$ is saturated. Let $x\in X$.
Since each $W_i$ is open, for each $i\in I$ one can write $W_i=\tbigcup_{j\in J_i} U_j\times V_j$ with $U_j$ open in $X$ and $V_j$ open in $Y$. Then $W_i\cap \{x\}\times Y=\tbigcup_{j\in J_i, x\in U_j}\{x\}\times V_j= \{x\}\times\left( \tbigcup_{j\in J_i, x\in U_j} V_j\right)$, and so $A\cap \{x\}\times Y=\tbigcap_{i\in I} \{x\}\times\left( \tbigcup_{j\in J_i, x\in U_j} V_j\right)=\{x\}\times\left( \tbigcap_{i\in I}\tbigcup_{j\in J_i, x\in U_j} V_j\right)$. Thus $\pi_2[A\cap\{x\}\times Y]=\tbigcap_{i\in I}\tbigcup_{j\in J_i, x\in U_j} V_j$ is fitted.
\end{proof}

\begin{remark}
In view of the result above, compactness with respect to the saturation operator imposes no restriction on a topological space. We warn the reader that this does \emph{not} imply that every spatial locale is $\FF$-compact. Although saturation and fitting are both defined by intersecting open subobjects, fitted sublocales form a non-conservative extension of saturated subspaces:  a spatial locale typically has nonspatial fitted sublocales.  Moreover, localic products of spatial locales in general do not coincide with the locales induced by the corresponding topological products, and the definition of $\FF$-compactness requires testing against arbitrary locales, including nonspatial ones. As we shall see, this makes $\FF$-compactness a strictly stronger condition, even for spatial locales.
\end{remark}

\section{The locale of pairwise disjoint sequences and semi-irreducibility}\label{sect4}
Semi-irreducible topological spaces were introduced by A.\,H. Stone in \cite{Stone1960} in his study of hereditarily compact spaces, where a space is \emph{semi-irreducible} if every pairwise disjoint family of (nonempty) open sets is finite.
Under classical logic, Dube extended the notion conservatively to point-free topology in \cite{irred}, where a locale is  \emph{semi-irreducible} if each of its pairwise disjoint subsets is finite.  With the aim of establishing a constructive theory of semi-irreducible locales, we now introduce one of the main ingredients of the paper:

\begin{definition}\label{def.Delta}
Let $L$ be an overt locale. The \emph{locale of pairwise disjoint sequences of positive opens of $L$} is the locale $\Delta_L$ presented by a generator $s_{n,a}$ for each $n\in\N$ and $a\in L^+$ subject to the following relations
\begin{enumerate}[label=\textup{(R\arabic*)},leftmargin=2.0\parindent]
\item\label{R1} $\tbigvee_{a\in L^+} s_{n,a}=1$ for each $n\in\N$;
\item \label{R2} $s_{n,a}\wedge s_{n,b}\leq \tbigvee\{1\mid a=b\}$ for each $n\in\N$ and $a,b\in L^+$;
\item \label{R3} $s_{n,a}\wedge s_{m,b}=0$ for every $n\ne m$ and $a\wedge b > 0$. 
\end{enumerate}
\end{definition}

\begin{remarks}
\begin{enumerate}[label=\textup{(\arabic*)},leftmargin=2.0\parindent]
\item  The points of the locale $\Delta_L$ are easily seen to be given by functions $f\colon \N\to L^+$ such that $f(n)\wedge f(m)>0$ is false for every $n\ne m$. By Lemma~\ref{property.overt} this holds if and only if $f(n)\wedge f(m)=0$ for every $n\ne m$---i.e., points correspond to pairwise disjoint sequences of positive opens of $L$. Indeed, each proposition $s_{n,a}$ is interpreted as $f(n)=a$, where \ref{R2} ensures that the relation is functional and \ref{R1}  ensures where the relation is total, whereas \ref{R3} provides pairwise disjointness. 
\item Classically, we shall show in Proposition~\ref{clasically.spatial} that the locale $\Delta_L$ is always spatial; hence  asking that $\Delta_L=\mathsf{O}$ is equivalent to asking that $L$ has no infinite pairwise disjoint sequences of nonzero elements (assuming classical logic). However, constructively, asking that $\Delta_L=\mathsf{O}$ is in general a stronger condition.
\end{enumerate}
\end{remarks}

\begin{definition}
A locale $L$ is \emph{geometrically semi-irreducible} if the locale $\Delta_L$ of pairwise disjoint sequences of positive opens of $L$ is trivial---i.e., if $\Delta_L=\mathsf{O}$.
\end{definition}

\begin{remarks}
\begin{enumerate}[label=\textup{(\arabic*)},leftmargin=2.0\parindent]
\item The reason for using the term `geometric' refers the fact that $\Delta_L$ classifies the propositional geometric theory of pairwise disjoint sequences. A different, but related, reason is motivated by Theorem~\ref{theorem.geometrically} below.
\item There are other sensible notions of semi-irreducibility constructively, that we shall not develop in this paper. Here, we just restrict to irreducibility in line with the notion of connected element of a locale (cf. \cite[Lemma~C.1.5.7]{elephant}), by considering a locale  \emph{irreducible} if for every pairwise disjoint family $\{a_i\}_{i\in I}$ of $L$ with $\tbigvee_{i\in I}a_i >0$, there is an $i\in  I$ such that $a_i=\tbigvee_{j\in I} a_j$. Here \emph{pairwise disjointness} is understood constructively as $a_i\wedge a_j\leq \tbigvee\{1\mid i=j\}$ for every $i,j\in I$ (see also \cite[pp.~1891--1892]{Vickers2009}). In this sense, a locale is irreducible if and only if each of its positive elements is connected.
\end{enumerate}
\end{remarks}

\section{Geometric semi-irreducibility implies $\FF$-compactness}\label{sect.imp1}

Our characterisation of $\FF$-compactness applies to overt locales. We do not know whether it extends to arbitrary locales or whether $\FF$-compactness itself implies overtness. With regards to  the second question, recall that a locale $L$ is overt precisely when the unique map from $L$ into the terminal locale is open. By contrast, it is not difficult to check that the unique map from $L$ into the terminal locale is  fitted for every locale $L$.

\begin{lemma}\label{c.BL.1}
Let $L$ be an overt locale, and $B$ a complete Boolean algebra. Then the composite $$q\colon \B_{L\oplus B}\hookrightarrow L\oplus B\twoheadrightarrow B$$ is open. Moreover, $q_! ( (a\oplus 1)^{**})=1$ for any $a>0$ in $L$.
\end{lemma}

\begin{proof}
The frame homomorphism corresponding to the map in the statement is $q^*:=(-)^{**}\circ \iota_2$. Since $\iota_2(b)$ is complemented for any $b\in B$, it is easy to check that  $q^*$ has a left adjoint given by $q_{!}(D)= \pi_{2!}(D)$ for any $D\in \B_{L\oplus B}$. Moreover, \eqref{fr} is clearly satisfied. The second part is a consequence of Lemma~\ref{descript.lo.pi}.
\end{proof}

\begin{lemma}\label{c.BL.2}
Let $L$ be an overt geometrically semi-irreducible locale,  $B$ a complete Boolean algebra  and denote $q$ as in Lemma~\ref{c.BL.1}. Let  $\{D_i\}_{i\in I}\subseteq \B_{L\oplus B}$. If $q_{!} (\tbigwedge_{i\in F} D_i)=1$ for any $F\in \K(I)$, then $q_!(\tbigwedge_{i\in I} D_i)=1$. 
\end{lemma}

\begin{proof}
Let $d:=q_!(\tbigwedge_{i\in I} D_i)^*$. Since $B$ is Boolean, our goal is to show $d=0$.  For each $n\in \N$ and $(a_1,\dots,a_n)\in (L^{+})^n$, let
$$d_{(a_1,\dots,a_n)}:= d\wedge \tbigwedge_{j=1}^{n} \tbigvee_{i\in I} q_{!} ((a_j\oplus 1)^{**}\wedge D_i)^* \wedge \tbigwedge_{j=1}^n\tbigwedge_{k=1}^{j-1} q_!( ((a_j\wedge a_k)\oplus 1)^{**})^*\in B.$$
Define 
$$E_0:= \op_B(d)\quad\text{ and } \quad E_n=\tbigcoprod^\cat{Loc}_{(a_1,\dots,a_n)\in (L^{+})^n}\op_B( d_{(a_1,\dots,a_n)})\quad\text{ for }n\in\N.$$
Note that since $d_{(a_1,\dots,a_{n},a_{n+1})}\leq d_{(a_1,\dots,a_n)}$ each $\op_B( d_{(a_1,\dots,a_{n+1})})$ embeds as an (open) sublocale of $\op_B( d_{(a_1,\dots,a_n)})$, and hence of $E_n$. We  thus get a  localic map $f_n\colon E_{n+1}\to E_n$ for each $n\geq 0$.  
Next we show that 
\begin{equation}\label{aim1}
d_{(a_1,\dots,a_n)}=\tbigvee_{a\in L^+} d_{(a_1,\dots,a_n,a)}\end{equation} for any $a_1,\dots,a_n\in L^+$.  Fix $a_1,\dots,a_n\in L^+$.

We introduce the following notation for $i,i_1,\dots i_n\in I$ and $j\leq n$:  
$$t_{i,j}:= q_{!}\left(  (a_j\oplus 1)^{**}\wedge D_i  \right)^*\quad\text{and}\quad t_{(i_1,\dots i_n)}:= q^*\left( d_{(a_1,\dots,a_n)}\wedge \tbigwedge_{k=1}^{n} t_{i_k,k}\right)\wedge \tbigwedge_{k=1}^{n} D_{i_k}.$$

We first show that  for all $i\in I$ and $(i_1,\dots,i_n)\in I^n$ the relation 
\begin{equation}\label{aim11}D_{i}^*\wedge t_{(i_1,\dots,i_n)}\leq \tbigvee_{a\in L^+} q^* (d_{(a_1,\dots,a_n,a)})\end{equation} holds. Since a basis of $\mathcal{B}_{L\oplus B}$ is given by elements of the form 
$$x_{ab}:=(a\oplus 1)^{**}\wedge q^*(b)$$ 
for $a\in L^+$ (by overtness of $L$) and $b\in B$, assume that $x_{ab}\leq D_{i}^*\wedge t_{(i_1,\dots,i_n)}$. In order to show \eqref{aim11}, it suffices to show $x_{ab}\leq q^*( d_{(a_1,\dots,a_n,a)})$, or by adjunction $q_!(x_{ab})\leq  d_{(a_1,\dots,a_n,a)}$.  Note that by \eqref{fr} and the second part of Lemma~\ref{c.BL.1}, one has $q_!(x_{ab})=b$. 
Moreover, $q_!(t_{(i_1,\dots,i_n)})= d_{(a_1,\dots,a_n)}\wedge \tbigwedge_{k=1}^{n} t_{i_k,k}$ using the assumption and \eqref{fr}; and so, $b=q_!(x_{ab})\leq q_!(t_{(i_1,\dots i_n)})\leq d_{(a_1,\dots,a_n)}$. 


Hence, in order to show that $b\leq d_{(a_1,\dots,a_n,a)}$, it suffices to show that $b\leq \tbigvee_{i\in I} q_!((a\oplus 1)^{**}\wedge D_i)^*$ and $b\leq q_!( ((a\wedge a_j)\oplus 1)^{**})^*$ for any $j\in\{1,\dots,n\}$. For the former, since $x_{ab}\wedge D_i=0$, by \eqref{fr}, it follows that $0=b\wedge q_!( (a\oplus 1)^{**}\wedge D_i)$ and hence $b\leq q_!( (a\oplus 1)^{**}\wedge D_i)^{*}$. For the latter, for any $j\in\{1,\dots,n\}$, by \eqref{fr},
$$b\wedge q_!( ((a\wedge a_j)\oplus 1)^{**}) = q_!(x_{ab}\wedge (a_j\oplus 1)^{**})\leq q_!(t_{(i_1,\dots,i_n)}\wedge (a_j\oplus 1)^{**})=0$$
because  \begin{align*}t_{(i_1,\dots , i_n)}\wedge (a_j\oplus 1)^{**}&\leq q^*(d_{(a_1,\dots,a_n)}\wedge t_{i_j,j})\wedge D_{i_j}\wedge (a_j\oplus 1)^{**}
\\& \leq q^*(d_{(a_1,\dots,a_n)}\wedge t_{i_j,j})\wedge q^*(q_!( D_{i_j}\wedge (a_j\oplus 1)^{**}))
\\&= q^* ( d_{(a_1,\dots,a_n)}\wedge t_{i_j,j}  \wedge q_!( D_{i_j}\wedge (a_j\oplus 1)^{**}))
\\& = q^*(d_{(a_1,\dots,a_n)}\wedge t_{i_j,j}\wedge  t_{i_j,j}^*)=0.
\end{align*}
Hence \eqref{aim11} follows.  Now, by adjunction, $t_{(i_1,\dots,i_n)}\wedge\tbigwedge_{i\in I}D_i\leq  q^*(q_! (t_{(i_1,\dots,i_n)}))\wedge q^*(q_!(\tbigwedge_{i\in I} D_i))=0$ because
$q_!(t_{(i_1,\dots,i_n)})\leq d_{(a_1,\dots,a_n)}\leq d=q_!(\tbigwedge_{i\in I}D_i)^*$. 
By Booleanness of the $\B_{L\oplus B}$ it follows that $t_{(i_1,\dots,i_n)}=\tbigvee_{i\in I}D_{i}^*\wedge t_{(i_1,\dots,i_n)}$. Combining it with \eqref{aim11}, it follows that 
$$t_{(i_1,\dots,i_n)}\leq \tbigvee_{a\in L^+} q^* (d_{(a_1,\dots,a_n,a)}),$$
and by adjunction, and applying \eqref{fr} and the hypothesis,
\begin{equation}d_{(a_1,\dots,a_n)}\wedge \tbigwedge_{k=1}^{n}t_{i_k,k}\leq \tbigvee_{a\in L^+}d_{(a_1,\dots,a_n,a)}\end{equation}
for any $i_1,\dots,i_n\in I$.
Finally,  using the distributivity of Kuratowski-finite meets over arbitrary joins one has
$ d_{(a_1,\dots,a_n)}\leq \tbigwedge_{k=1}^n \tbigvee_{i\in I}t_{i,k}= \tbigvee_{(i_1,\dots,i_n)\in I^n}\tbigwedge_{k=1}^{n}t_{i_k,k}$, from which \eqref{aim1} follows.  

By Lemma~\ref{lemma.open.coproducts} and \eqref{aim1}, it is easy to check that localic map $\tbigcoprod^{\cat{Loc}}_{a\in L^+}\op_B(d_{(a_1,\dots,a_n,a)}) \to \op_B(d_{(a_1,\dots,a_n)})$ is surjective, and therefore the localic maps $f_n\colon E_{n+1}\to E_n$ are surjective for all $n\geq 0$. 
By \cite[Corollary~C.1.1.12]{elephant} it follows that the localic maps $p_n\colon \varprojlim_{n\geq 0} E_n \to E_n$ are surjective for each $n\geq 0$.

For each $a\in L^+$ and $n\in \N$, let 
$$u_{n,a}=  ( d_{(a_1,\dots,a_n)}\to \tbigvee^{B}\{ 1\mid a_n=a\})_{(a_1\dots,a_n)\in (L^+)^n}\in E_n.$$

We define a frame homomorphism $h\colon \Delta_L\to  \varprojlim_{n\geq 0} E_n$ given by $h(s_{n,a})= p_{n}^*(u_{n,a})$ for $n\in \N$ and $a\in L^+$. We check that $h$ sends the relations \ref{R1}-\ref{R3} into identities in $ \varprojlim_{n\geq 0} E_n$. 

For $\ref{R1}$, it suffices to show that $\tbigvee_{a\in L^+}u_{n,a}=1$. Since joins in the cartesian product are given componentwise, one must show that for any $(a_1,\dots,a_n)\in (L^+)^n$, the relation $\tbigvee_{a\in L^+}^{\op_B(d_{(a_1,\dots,a_n)})}(d_{(a_1,\dots,a_n)}\to \tbigvee^B\{1\mid a_n=a\})=1$, which clearly holds because 
$\tbigvee^{\op_B(d_{(a_1,\dots,a_n)})}_{a\in L^+}(d_{(a_1,\dots,a_n)}\to \tbigvee^B\{1\mid a_n=a\})\geq d_{(a_1,\dots,a_n)}\to \tbigvee^B\{1\mid a_n=a_n\}= d_{(a_1,\dots,a_n)}\to 1=1$.

Let us now show $\ref{R2}$ so let $a,b\in L^+$ and $n\in\N$. It suffices to show that $u_{n,a}\wedge u_{n,b}\leq \tbigvee^{E_n}\{ 1\mid a=b\}$, and since joins and meets are given componentwise, it is equivalent to show that for any $(a_1,\dots,a_n)\in (L^+)^n$, the relation $d_{(a_1,\dots,a_n)}\to \tbigvee^{B}\{ 1\mid a_n=a\}\wedge d_{(a_1,\dots,a_n)}\to \tbigvee^{B}\{ 1\mid a_n=b\}\leq  \tbigvee^{\op(d_{(a_1,\dots,a_n))}}\{ 1\mid a=b\}$; but this clearly holds because the surjection $B\twoheadrightarrow \op_B(d_{(a_1,\dots,a_n)})$ is a frame homomorphism.

Finally we show $\ref{R3}$. Let $n<m$ and $a,b\in L^+$ with $a\wedge b>0$. Consider the diagram 
\[\begin{tikzcd}
	{E_n} & {E_m} \\
	{\varprojlim_{n\geq 0} E_n}
	\arrow["{f_{nm}^*}", from=1-1, to=1-2]
	\arrow["{p_{n}^*}"', from=1-1, to=2-1]
	\arrow["{p_{m}^*}", from=1-2, to=2-1]
\end{tikzcd}\]
It therefore suffices to show that $f_{nm}^*(u_{n,a})\wedge u_{m,b}=0_{E_m}$.  Moreover, it is easy to check that the transition frame homomorphism $f_{nm}^*$ sends $u_{n,a}$ to the element $(  d_{(a_1\dots,a_m)} \to \tbigvee^B \{1\mid a_n=a\})_{(a_1,\dots,a_m)\in (L^+)^m}$. Hence, we are left with the task of showing that for any $(a_1,\dots,a_m)\in (L^+)^m$, the relation $$(d_{(a_1\dots,a_m)} \to \tbigvee^B \{1\mid a_n=a \})\wedge (d_{(a_1\dots,a_m)} \to \tbigvee^B \{1\mid a_m=b\})=0_{\op_B(d_{(a_1,\dots,a_m)})}$$
or equivalently, that 
$$\tbigvee^{\op_B(d_{(a_1,\dots,a_m)})} \{ 1\mid a_n=a,\ a_m=b\}\leq  0_{\op_B(d_{(a_1,\dots,a_m)})}$$
---i.e., that if $a_n=a$ and $a_m=b$ then $d_{(a_1,\dots,a_m)}=0$.  But this is clear because $d_{(a_1,\dots,a_m)}\leq q_! (( (a\wedge b)\oplus 1)^{**})^*$ and because of Lemma~\ref{c.BL.1} and the fact that $a\wedge b>0$. Hence $h$ is a frame homomorphism. Since $1=0$ in $\Delta_L$, it follows that $1=0$ in $\varprojlim_{n\geq 0} E_n$--- i.e. $\varprojlim_{n\geq 0} E_n=\mathsf{O}$. Since $p_0\colon \varprojlim_{n\geq 0} E_n\to E_0$ is a surjection, it follows that $E_0=\mathsf{O}$, and hence $d=0$, that is, $q_!(\tbigwedge_{i\in I}D_i)=1$ by Booleanness of $B$.
\end{proof}

\begin{theorem}\label{mainfir}
Let $L$ be an overt geometrically semi-irreducible locale. Then $L$ is $\FF$-compact.
\end{theorem}

\begin{proof}
Let $L$ be overt, $M$ arbitrary and $\{D_i\}_{i\in I}\subseteq L\oplus M$. We will show that $$\pi_2[ \tbigcap_{i\in I}\op_{L\oplus M}(D_i)]= \tbigcap_{F\in \K(I)} \op_M(\pi_{2!} (\tbigwedge_{i\in F}D_i)).$$
The inclusion $\subseteq$ is trivially satisfied because of the adjunction $\pi_{2!}\dashv \pi_{2}^{*}$, so we prove the reverse inclusion. 
Denote  $R:=\tbigcap_{F\in \K(I)} \op_M(\pi_{2!} (\tbigwedge_{i\in F}D_i))$ and $\iota\colon R\hookrightarrow M$ the sublocale inclusion. Moreover, let $$p_R\colon \SS_b(R)\twoheadrightarrow  R$$ the localic surjection from Lemma~\ref{sbl}. For any $F\in \K(I)$, the map $ip_R$ factors through $\op_M(\pi_{2!}(\tbigwedge_{i\in F}D_i))$ from which follows from Lemma~\ref{lemma.factors} that $(ip_R)^*(  \pi_{2!}(\tbigwedge_{i\in F}D_i) )= 1$. 
Consider the pullback diagram
\[\begin{tikzcd}
	{L\oplus \SS_b(R)} & {\SS_b(R)} \\
	{L\oplus M} & M
	\arrow["\pi"', from=1-1, to=1-2]
	\arrow["{h:=1_L\oplus \iota p_R}"', from=1-1, to=2-1]
	\arrow["{\iota p_R}"', from=1-2, to=2-2]
	\arrow["{\pi_2}", from=2-1, to=2-2]
\end{tikzcd}\]
Since $\pi_2$ is open by overtness (cf. Proposition~\ref{overt.projec}), the Beck-Chevalley condition $$(\iota p_R)^*\circ \pi_{2!}=\pi_{!} \circ h^*$$ holds
(cf. \cite[Proposition~C3.1.11]{elephant}). Hence, for each $F\in\mathcal{K}(I)$, one has $\pi_! (\tbigwedge_{i\in F} h^* (D_i ))= \pi_! (h^* (\tbigwedge_{i\in F}D_i ))=1$.  Letting $D_{i}'= h^*(D_i)^{**}\in \mathcal{B}_{L\oplus \SS_b(R)}$, by Lemma~\ref{c.BL.1} it follows that $q_!(\tbigwedge_{i\in F} D_{i}')=1$ where $q\colon  \mathcal{B}_{L\oplus \SS_b(R)}\to \SS_b(R)$ is the projection. By Lemma~\ref{c.BL.2}, one has $q_!(\tbigwedge_{i\in I}D_{i}^{'})=1$.  Note that the composite
\[\begin{tikzcd}
	{g\colon \op_{  \mathcal{B}_{L\oplus \SS_b(R)} }( \tbigwedge_{i\in I} D_{i}^{'} )} & { \mathcal{B}_{L\oplus \SS_b(R)}} & {\SS_b(R)}
	\arrow["j", from=1-1, to=1-2]
	\arrow["q", from=1-2, to=1-3]
\end{tikzcd}\]is open by Lemma~\ref{c.BL.1} and since $g_!(1)= q_!j_!(1)=q_! ( \tbigwedge_{i\in I} D_{i}^{'}  )=1$; it follows from Lemma~\ref{open.surj} that $g$ is surjective. Finally, if we denote by $j'\colon  \mathcal{B}_{L\oplus \SS_b(R)} \hookrightarrow L\oplus \SS_b(R)$ the inclusion, then $(hj'j)^*(D_i)=j^*(D_{i}')=1$ for any $i\in I$,  so by Lemma~\ref{lemma.factors} it follows that $\text{Im}(hj'j)\subseteq \tbigcap_{i\in I}\op_{L\oplus M} (D_i)$ and so $\text{Im}(\pi_2hj'j)\subseteq \pi_2[\tbigcap_{i\in I}\op_{L\oplus M} (D_i)]$. On the other hand $\pi_2 h j' j= \iota  p_R \pi  j' j=\iota p_R  q  j = \iota p_R g$, and since $g$ and $p_R$ are surjective,  it follows that $\text{Im} (\pi_2hj'j)=R$, which establishes the required inclusion.\end{proof}

\section{The converse implication}\label{sect.imp2}

\begin{lemma}\label{conv1}
Let $B$ be a complete Boolean algebra. Then there exists a complete Boolean algebra $ \hat{B}$ and an injective frame homomorphism $\eta\colon B\to \hat{B}$ together with elements $a_{n,k}\in  \hat{B}$ for each $(n,k)\in\N\times\N$ such that 
\begin{enumerate}[label=\textup{(\arabic*)},leftmargin=2.0\parindent]\item\label{cantor0} $B=\mathsf{O}$ if and only if $\hat B= \mathsf{O}$;
\item \label{cantor1} $a_{n,1}=1$, $a_{n,k+1}\leq a_{n,k}$ and $\bigwedge_{k\in\N}a_{n,k}=0$ for each $n,k\in \N$;
\item  \label{cantor2} for each $n\in\N$ the map $\lambda_n\colon  \hat{B}\to B$ given by $$\lambda_n(c)=\tbigvee\{b\in B\mid \exists k\in\N\text{ such that } \eta(b)\wedge a_{n,k}\leq c\}$$ preserves Kuratowski-finite meets;
\item \label{cantoradded2} If $b\in B$, then $b\leq\lambda_n(\eta(b))$ for any $n\in\N$;
\item  \label{cantor3} If for all $n\in\N$ the relation $\lambda_n(c)=1$ holds, then $c=1$.
\end{enumerate}
\end{lemma}

\begin{proof}
 Let $K$ be the Cantor frame on the set $\N\times\N$, namely the frame presented by symbols $z_{n,k}$ and $t_{n,k}$ for each $(n,k)\in\N\times\N$ subject to the relations (cf. \cite[Example~2.2]{Manuell2026}):
\begin{enumerate}[label=\textup{(K\arabic*)},leftmargin=2.0\parindent]
\item $z_{n,k}\vee t_{n,k}=1$ for each $n,k\in\N$,
\item $z_{n,k}\wedge t_{n,k}=0$ for each $n,k\in\N$.
\end{enumerate}
 The frame $ \hat{B}:=\mathcal{B}_{B\oplus K}$ and the elements $a_{n,k}:= \tbigwedge_{m<k} 1\oplus z_{n,m}\in  \hat{B}$ satisfy the required conditions. Note that the composite \[\begin{tikzcd}
	{\eta:=B} & {B\oplus K} & {\mathcal{B}_{B\oplus K}}
	\arrow["{\iota_1}", hook, from=1-1, to=1-2]
	\arrow["q", two heads, from=1-2, to=1-3]
\end{tikzcd}\]
is  injective, because $B$ is Boolean and $\iota_1$ preserves complements. Moreover, \ref{cantor0} is trivial and for \ref{cantor1}, the first two identities are also obvious, so let us show that 
$\bigwedge_{k\in\N}a_{n,k}=0$ for each $n,k\in \N$. A basis for $K$ is given by Kuratowski-finite meets of generators---i.e., by elements of the form $\tbigwedge F$ where 
$$F\in \K(\{z_{n,k}\mid n,k\in\N\}\cup \{t_{n,k}\mid n,k\in\N \}).$$
 Any such $F$ can be decomposed as $F=P\cup Q$ where 
 $$P=\{ x\in F\mid \exists n,k\in\N,\ x={z_{n,k}}\}\quad\text{and}\quad Q=\{ x\in F\mid \exists n,k\in\N,\ x={t_{n,k}}\}.$$ 
 Since $P\cap Q=\varnothing$, both subsets are complemented and hence Kuratowski-finite by Lemma~\ref{kf}\,\ref{kf1}. Thus, an alternative basis is given by elements of the form $\tbigwedge_{(n,k)\in P} z_{n,k}\wedge \tbigwedge_{(n,k)\in Q} t_{n,k}$ for Kuratowski-finite subsets $P,Q\subseteq \N\times \N$. 
Hence, a basis for $ \hat{B}$ is given by elements of the form 
$$x_{b,P,Q}:=\iota_1(b)\wedge \iota_2( \tbigwedge_{(n,k)\in P} z_{n,k}\wedge \tbigwedge_{(n,k)\in Q} t_{n,k})$$ for $P,Q\in \K(\N\times \N)$. Hence, let any such basic element $x_{b,P,Q}$ below $\tbigwedge_{k\in\N} \tbigwedge_{m<k}1\oplus z_{n,m}$ for some $n\in\N$. Since $P$ and $Q$ are Kuratowski-finite, choose $(n_0,m_0)\in\N\times\N$, such that for all $(n',m')\in P$ and $(n'',m'')\in Q$, one has $n',n''<n_0$ and $m',m''<m_0$ (the set $P\cup Q$ is Kuratowski-finite; if it is empty, choose $(n_0,m_0)=(1,1)$, otherwise take the coordinatewise sum and add $1$, cf. Lemma~\ref{kf}\,\ref{kf2}). Then $x_{b,P,Q}\leq 1\oplus z_{n,m_0}$ and so $x_{b,P,Q}\wedge 1\oplus t_{n,m_0}=0$.  Consider the frame homomorphism $f\colon K\to K$ such that $f(z_{n,m_0})=0$ and $f(t_{n,m_0})=1$ and fixing the remaining generators. Applying the frame homomorphism $1_B\oplus f$, it then follows that $x_{b,P,Q}=0$. Hence $\bigwedge_{k\in\N}a_{n,k}=0$.\\[2mm]
\ref{cantor2}. Clearly, $\lambda_n$ preserves the top element. Assume now that there are $k,k'\in\N$ such that $\eta(b)\wedge a_{n,k} \leq c$ and $\eta(b')\wedge a_{n,k'}\leq d$. Then, $\eta (b\wedge b')\wedge a_{n,\text{max}\{k,k'\}}\leq c\wedge c'$. Hence $\lambda_{n}(c\wedge c')\leq  \lambda_n(c)\wedge\lambda_n(c')$, as required.\\[2mm]
\ref{cantoradded2} is straightforward from the definition of $\lambda_n$ (choose $k=0$).\\[2mm]
\ref{cantor3}  Let $c\in  \hat{B}$ with $\lambda_n(c)=1$ for all $n\in\N$. Let $x_{b,P,Q}$ be a basic generator with $P,Q\in\mathcal{K}(\N\times \N)$  such that $x_{b,P,Q}\leq c^*$.  If we show that $x_{b,P,Q}=0$, we are done by Booleanness of $\hat B$. Since $P$ and $Q$ are finite, choose $(n_0,m_0)\in\N\times\N$, such that for all $(n',m')\in P$ and $(n'',m'')\in Q$, one has $n',n''<n_0$ and $m',m''<m_0$. Let $b'\in B$ and $k\in\N$ such that $\eta(b')\wedge a_{n_0,k}\leq c$.  Then, $(b'\oplus 1)\wedge \tbigwedge_{m<k} 1\oplus z_{n_0,m} \wedge x_{b,P,Q}=0$. 
Consider the frame homomorphism $h\colon K\to K$ given by $h(z_{n_0,k'})=1$ and $h(t_{n_0,k'})=0$ for all $k'\in K$ and which fixes the rest of generators. Applying $1_B\oplus h$, we get $b'\oplus 1\wedge x_{b,P,Q}=0$ and so $\lambda_{n_0}(c)\oplus 1 \wedge x_{b,P,Q}=0$. Since $\lambda_{n_0}(c)=1$, the conclusion follows.
\end{proof}

In what follows we recall the \emph{Artin glueing} for frames (see, for instance, \cite{FaulManuell2020} or \cite[Section~A.4.5]{elephant}). The Artin glueing of frames $H$ and $N$ along a Kuratowski-finite-meet preserving map $\lambda\colon H\to N$  is given by the frame
$$\text{Gl}(\lambda):= \{ (n,h)\in N\times H\mid n\leq \lambda(h) \}$$
with arbitrary joins and Kuratowski-finite meets given componentwise. There are projection frame surjections $p_H\colon \text{Gl}(\lambda)\to H$, corresponding to an open sublocale, and $p_N\colon \text{Gl}(\lambda)\to N$, corresponding to its closed complement.

The following are easy properties of the Artin glueing construction:

\begin{lemma}\label{basic.artin}
Let $H$ and $N$ be frames and $\lambda\colon H\to N$ a meet-preserving map. Then,
\begin{enumerate}[label=\textup{(\arabic*)},leftmargin=2.0\parindent]
\item\label{basic.artin1} The Heyting operation in $\text{Gl}(\lambda)$ is given by $$(n,h)\to (n',h')= ( (n\to n')\wedge \lambda (h\to h'), h\to h')$$ for every $(n,h),(n',h')\in\text{Gl}(\lambda)$;
\item\label{basic.artin2} The right adjoint of the projection $p_N$ is given by $p_{N*}(n)=(n,1)$ for every $n\in N$.
\end{enumerate}
\end{lemma}

\begin{lemma}\label{artinglueing}
Let $B$ be a complete Boolean algebra and let $ \hat{B}$  be the complete Boolean algebra and $\lambda_n\colon  \hat{B}\to B$ the finite-meet preserving maps   satisfying the conditions of Lemma~\ref{conv1}. Then the map $\lambda\colon  \hat{B}\to \N{\cdot} B$ given by $\lambda(c)=(\lambda_n(c))_{n\in\N}$ preserves Kuratowski-finite meets, and the Artin glueing $\text{Gl}(\lambda)$ of $ \hat{B}$ and $ \N{\cdot} B$ satisfies the following properties:
\begin{enumerate}[label=\textup{(\arabic*)},leftmargin=2.0\parindent]\item\label{artin1} $p_{ \N{\cdot} B*}[ \N{\cdot} B]^\circ = \text{Gl}(\lambda)$;
\item\label{artin2}   $p_{ \N{\cdot} B*}[ \N{\cdot} B]$ is fitted in $\text{Gl}(\lambda)$ if and only if $B=\mathsf{O}$;
\item\label{artin3}  If $\iota_n\colon B\to  \N{\cdot} B$ denotes the localic injection into the coproduct, then $p_{ \N{\cdot} B*}\iota_{n}[B]$ is fitted in  $\text{Gl}(\lambda)$ for any $n\in\N$;
\item\label{artin4} The map $p^*\colon B\to \text{Gl}(\lambda)$ given by $p^*(b):= ( (b)_{n\in\N}, \eta(b) )$ is an injective frame homomorphism.
\end{enumerate}
\end{lemma}

\begin{proof}
Since frame operations in $ \N{\cdot} B$ are given componentwise, it is clear that $\lambda$ preserves finite meets because do does each $\lambda_n$.\\[2mm]
\ref{artin1} Let $((b_n)_{n\in\N},c)\in \text{Gl}(\lambda)$ such that   $p_{ \N{\cdot} B*}[ \N{\cdot} B]\subseteq\op( (b_n)_{n\in\N},c)$, and so $b_n=1$ for all $n\in\N$.  Since $((b_n)_{n\in\N}),c)\in\text{Gl}(\lambda)$, it follows that $\lambda_n(c)=1$ for all $n\in \N$. From Lemma~\ref{conv1}, one concludes $c=1$, as required.\\[2mm]
\ref{artin2} By virtue of \ref{artin1},  $p_{ \N{\cdot} B*}[ \N{\cdot} B]$ is fitted if and only if $p_{ \N{\cdot} B*}[ \N{\cdot} B]=\text{Gl}(\lambda)$ and since $ \hat{B}$ is its open complement, the latter holds iff $ \hat{B}=\mathsf{O}$. By Lemma~\ref{conv1}\,\ref{cantor0} it is clear that the last condition holds if and only if $B=\mathsf{O}$.\\[2mm]
\ref{artin3} Fix $n\in\N$ and denote by $(b_m)_{m\in\N}\in  \N{\cdot} B$ the element given by $b_n=1$ and $b_{m}=0$ for $m\ne n$ (note that constructively it is well-defined because $\N$ has decidable equality). Note that for any $k\in\N$ one has  $( (b_m)_{m\in\N} , a_{n,k})\in\text{Gl}(\lambda)$ because $\lambda_n(a_{n,k})=1$.  We claim that $p_{ \N{\cdot} B*}\iota_{n}[B]=\tbigcap_{k\in\N}\op( ((b_m)_{m\in\N} , a_{n,k}))$. The inclusion $\subseteq$ follows easily by adjunction so let us prove the reverse one. Let $( (a_m)_{m\in\N},c)\in \tbigcap_{k\in\N}\op( ((b_m)_{m\in\N} , a_{n,k}))$. By Lemma~\ref{basic.artin}\,\ref{basic.artin1}, this means  $a_{n,k}\to c=c$  and $(b_m\to a_m)\wedge \lambda_m(c)\leq a_m$ for every $k,m\in\N$.  From the first equality, one has $a_{n,k}^*\leq c$ for every $k\in \N$ and so $\tbigvee_{k\in\N} a_{n,k}^*\leq c$. By Booleanness of $\hat B$ and Lemma~\ref{conv1}\,\ref{cantor1}, it follows that $c=1$. The second equality then implies that $b_m\to a_m\leq a_m$ for every $m\in \N$ and by  definition of $(b_m)_{m\in\N}$ it follows that  $a_m=1$ for any $m\ne n$. Therefore, $( (a_m)_{m\in\N}, c)\in p_{ \N{\cdot} B*}\iota_n[B]$ by the description of the right adjoint in Lemma~\ref{basic.artin}\,\ref{basic.artin2}.\\[2mm]
\ref{artin4} Note that $p^*$ is well defined by Lemma~\ref{conv1}\,\ref{cantoradded2}. Injectivity is trivial and the fact that $p^*$ is a frame homomorphism follows because frame operations in $\text{Gl}(\lambda)$ are given componentwise and $\eta$ is a frame homomorphism. 
\end{proof}

\begin{theorem}\label{mainsec}
Let $L$ be an overt $\FF$-compact locale. Then $L$ is geometrically semi-irreducible. 
\end{theorem}

\begin{proof}
Let $B=\mathcal{B}_{\Delta_L}$, construct the Artin glueing $\text{Gl}(\lambda)$  of $\hat B$ and $ \N{\cdot} B$ as in Lemma~\ref{artinglueing}, and let $p\colon \text{Gl}(\lambda)\twoheadrightarrow B$ the localic  surjection corresponding to the injective frame homomorphism from  Lemma~\ref{artinglueing}\,\ref{artin4}.

Throughout the proof we denote $\iota_1\colon L\to L\oplus \text{Gl}(\lambda)$ and $\iota_2\colon \text{Gl}(\lambda)\to L\oplus \text{Gl}(\lambda)$ the coproduct injections and $\pi\colon L\oplus \text{Gl}(\lambda)\to \text{Gl}(\lambda)$ the localic projection.

For each $n\in\N$, set 
$$U_n:=\tbigvee_{a\in L^+} a\oplus p^*(s_{n,a}^{**} ) \in L\oplus \text{Gl}(\lambda).$$
First, note that $U_n\wedge U_m=\mathsf{O}$ for $n\ne m$. Indeed, we use overtness of $L$ as follows:
\begin{align*}
U_n\wedge U_m &=\tbigvee_{a,b\in L^+} \iota_1(a\wedge b)\wedge \iota_2(p^* ((s_{n,a}\wedge s_{m,b})^{**})) 
\\& =\tbigvee_{a,b\in L^+} \iota_1 (\tbigvee\{ a\wedge b \mid a\wedge b>0\}) \wedge  \iota_2(p^* ((s_{n,a}\wedge s_{m,b})^{**})) 
\\&= \tbigvee_{a,b\in L^+} \tbigvee_{a\wedge b>0} \iota_1 (a\wedge b) \wedge  \iota_2(p^* ((s_{n,a}\wedge s_{m,b})^{**}))=0 \end{align*}
because each element in the last join equals zero by \ref{R3}.
Let $$S_n:= \pi_{-1}[C_n]\cap\op(U_n),$$ where $C_n:=p_{ \N{\cdot} B*}\iota_{n}[B]$, which is fitted in $\text{Gl}(\lambda)$ by Lemma~\ref{artinglueing}\,\ref{artin3}. Hence $S_n$ is fitted in $L\oplus \text{Gl}(\lambda)$ for each $n\in\N$. 
By Lemma~\ref{fitted.disjoint}, the join $\tbigvee_{n\in\N}S_n$ is fitted in $L\oplus \text{Gl}(\lambda)$. Hence, since $L$ is $\FF$-compact it follows that $\pi[\tbigvee_{n\in\N}S_n]=\tbigvee_{n\in \N} \pi[S_n]$ is fitted in $\text{Gl}(\lambda)$. 
Note that the restriction \[\begin{tikzcd}
	{\op(U_n)} & {L\oplus \text{Gl}(\lambda)} & {\text{Gl}(\lambda)}
	\arrow["\iota",hook, from=1-1, to=1-2]
	\arrow["\pi", from=1-2, to=1-3]
\end{tikzcd}\]
is open by overtness of $L$ and hence surjective by Lemma~\ref{open.surj}  because
\begin{align*}\pi_!(\iota_!(1))&=\pi_!(U_n)=\tbigvee_{a\in L^+} \pi_! (  a\oplus ( p^*(s_{n,a}^{**} ))  )= \tbigvee_{a\in L^+} p^*(s_{n,a}^{**})=p^*(\tbigvee_{a\in L^+}s_{n,a}^{**})=1
\end{align*}
by \ref{R1} and Lemma~\ref{descript.lo.pi}. 
Hence, since open surjections are pullback stable, it follows that the left vertical arrow in the pullback diagram
\[\begin{tikzcd}
	{S_n=(\pi\iota)_{-1}[C_n]} & {\op(U_n)} \\
	{C_n} & {\text{Gl}(\lambda)}
	\arrow[hook, from=1-1, to=1-2]
	\arrow[two heads, from=1-1, to=2-1]
	\arrow["{\pi\iota}", two heads, from=1-2, to=2-2]
	\arrow[hook, from=2-1, to=2-2]
\end{tikzcd}\]
is surjective---i.e., $\pi[S_n]=C_n$. Thus, $\tbigvee_{n\in\N}C_n$ is fitted. By Lemma~\ref{artinglueing}\,\ref{artin2}, one has $B=\mathsf{O}$, and since $B$ is dense in $\Delta_L$, it follows that $\Delta_L=\mathsf{O}$. 
\end{proof}

Combining Theorem~\ref{mainfir} and Theorem~\ref{mainsec} we conclude the main result of the paper:

\begin{corollary}\label{MAINRESULT}
An overt locale is $\FF$-compact if and only if it is geometrically semi-irreducible.
\end{corollary}
\section{The discrete case}\label{section.discrete}

In what follows, we consider the notions  introduced above and the constructively valid results internally in toposes; we shall freely work in the internal logic of these toposes. If $L$ is an overt locale in a topos $\mathcal{E}$, for clarity, we write $\Delta^\mathcal{E}_{L}$ for the locale of pairwise disjoint sequences of positive opens of $L$  internal to the topos $\mathcal{E}$ (cf. Definition~\ref{def.Delta}). Similarly, we denote by $\N_\EE$ the natural numbers object in $\EE$. If $f\colon \FF \to \EE$ is a geometric morphism, we recall here that $f^*(\N_\EE)$ is a natural numbers object in $\FF$, and  hence it is isomorphic to $\N_\FF$ (cf. \cite[Lemma~A.4.1.14]{elephant}), and if $D$ is an internal lattice in $\EE$, the object $f^*(D)$ is also an internal lattice in $\FF$ with finite lattice operations obtained by applying $f^*$ to the corresponding  operations of $L$ (cf. \cite[Section~IV.8]{MacLaneMoerdijk}).

 We also recall (see \cite[Section~2.2.10]{Henry2016}) that if $L$ is a locale internal to $\EE$, there is a locale $f^\#(L)$ internal to $\FF$ making the following diagram a pullback of toposes:
 
 \[\begin{tikzcd}
	{\text{Sh}_\FF( f^\#(L))} & \FF \\
	{\text{Sh}_\EE(L)} & \EE
	\arrow[hook, from=1-1, to=1-2]
	\arrow[from=1-1, to=2-1]
	\arrow["f", from=1-2, to=2-2]
	\arrow[hook, from=2-1, to=2-2]
\end{tikzcd}\]
 
 In the particular case of the locale   of pairwise disjoint sequences,  the next lemma provides a presentation  of the locale obtained under change of base:

\begin{lemma}\label{lemma.present}
Let $L$ be an overt locale in a topos $\mathcal{E}$ and $f\colon \mathcal{F}\to\mathcal{E}$ a geometric morphism. Then, $f^\#(\Delta^\mathcal{E}_{L})$ is presented by generators $r_{n,a}$ for $n\in \N_\mathcal{F}$ and  $a\in f^*( L^+)$ subject to the relations
\begin{enumerate}[label=\textup{(T\arabic*)},leftmargin=2.0\parindent]
\item\label{T1} $\tbigvee_{a\in f^*(L^+)} r_{n,a}=1$ for each $n\in \N_\FF$;
\item \label{T2} $r_{n,a}\wedge r_{n,b}\leq \tbigvee\{1\mid a=b\}$ for each $n\in \N_\FF$ and $a,b\in f^*(L^+)$;
\item \label{T3} $r_{n,a}\wedge r_{m,b}=0$ for every $n\ne m$ and $a,b\in f^*(L^+)$ with $a\wedge b\in f^*(L^+)$. 
\end{enumerate}
\end{lemma}

\begin{proof}
It is an immediate consequence of \cite[Corollary~5.4]{Vickers2004}, for which we express  the presentation given by \ref{R1}--\ref{R3} in the notation from  \cite[p.~388]{Vickers2004}. We define the following objects in $\EE$:
$$G:= \N_\EE\times L^+,\quad \text{and}\quad R:= \N_\EE \sqcup ( \N_\EE\times L^+\times L^+)\sqcup ( \text{Diag}_{\N_\EE}^{c} \times  P_L)$$
where $\text{Diag}_{\N_\EE}^{c}$ stands for the complement of the diagonal (note that $\N_\EE$ has decidable equality---i.e., its diagonal is complemented) and $P_L$ is given by the following pullback in $\EE$:\[\begin{tikzcd}
	{P_L} & {L^+} \\
	{L^+\times L^+} & L
	\arrow[from=1-1, to=1-2]
	\arrow[from=1-1, to=2-1]
	\arrow[hook, from=1-2, to=2-2]
	\arrow["\wedge", from=2-1, to=2-2]
\end{tikzcd}\]
Moreover,
$$\lambda \colon R\to \K(G)$$ is given by
$\lambda(n)=\varnothing$ for $n\in\N_\EE$, $\lambda(n,a,b)=\{(n,a),(n,b)\}$ for $n\in\N_\EE$ and $a,b\in L^+$, and $\lambda(n,m,a,b)= \{ (n,a),(m,b)\}$ for all $n\ne m$ and $a\wedge b>0$. 
$$\pi\colon D:=\tbigcoprod_{r\in R} D_r\to R$$
is the projection, where  $D_n:= L^+$ for $n\in\N_\EE$, $D_{n,a,b}:=\{ *\mid a=b\}$ for $n\in \N_\EE$ and $a,b\in L^+$, and $D_{n,m,a,b}:=\varnothing$, for $n\ne m$ and $a\wedge b>0$.  Moreover, $$\rho\colon D\to \K(G)$$ is given by the coproduct of the maps $D_n\to \K(G)$ sending $a\in L^+$ to $\{ (n,a)\}$, $D_{n,a,b}\to \K(G)$ sending $*$ to $\varnothing$ and the unique map from $\varnothing$. Now, the presentation for $f^\#(\Delta_{L}^\EE)$ is obtained by applying $f^*$ to all the components $(R,G,D,\pi,\lambda,\rho)$ of the presentation. Since $f^*$ preserves finite colimits, finite limits (hence the diagonal and its complement too), natural numbers objects and Kuratowski-finite powersets, it follows that the presentation is given precisely by \ref{T1}--\ref{T3}.
\end{proof}

\begin{lemma}\label{first.lemma.geom}
Let $f\colon \FF\to\EE$ be a geometric morphism between toposes. If $X$ is an object of $\EE$ and there is a monomorphism $\iota\colon \N_\FF\rightarrowtail  f^*(X)$  in $\FF$, then the locale $f^\#(\Delta_{\mathcal{P}(X)}^\EE)$ has a point.
\end{lemma}

\begin{proof}
By Lemma~\ref{lemma.present}, in order to define a frame homomorphism $ f^\#(\Delta_{\mathcal{P}(X)})\to \Omega_\FF$, it suffices to define a morphism $$h\colon\{ r_{n,a}\mid (n,a)\in \N_\FF\times f^*(\mathcal{P}(X)^+)\}\to \Omega_\FF$$ which sends the relations \ref{T1}--\ref{T3} into identities in $\Omega_\FF$. Let $s_X\colon X\to \mathcal{P}(X)^+$ be the map given by $s_X(x):=\{x\}$ and set $$h(r_{n,a})=\llbracket f^*(s_X ) (\iota(n)) =   a\rrbracket.$$ 
Clearly, \ref{T1} and \ref{T2} are satisfied.  Let us now show \ref{T3}.  Denote by $g:=\wedge\circ s_X\times s_X\colon X\times X \to \mathcal{P}(X)$ the map given by $g(x,y)=\{x\}\cap\{y\}$. One has $g^{-1}[\mathcal{P}(X)^+]\subseteq \text{Diag}_X$, and since $f^*$ preserves pullbacks, it follows  $f^*(g)^{-1}[f^*(\mathcal{P}(X)^+)]\subseteq \text{Diag}_{f^*(X)}$. Hence, since $\iota$ is a monomorphism, one has $n= m$ if and only if $\iota(n)= \iota(m)$ if and only if $f^*(s_X)(\iota(n))\wedge f^*(s_X)(\iota(m))\in f^*(\mathcal{P}(X)^+)$. Then, let 
 $a,b\in f^*(\mathcal{P}(X)^+)$ such that  $a\wedge b\in f^*(\mathcal{P}(X)^+)$ and $n\ne m$
and suppose that $f^*(s_X ) (\iota(n)) =   a$ and $f^*(s_X ) (\iota(m)) =   b$. Therefore, $n=m$ and so $\llbracket f^*(s_X ) (\iota(n)) =   a\rrbracket\wedge \llbracket f^*(s_X ) (\iota(m)) =   a\rrbracket =\bot$. 
\end{proof}

\begin{lemma}\label{lemma.QQN}
Let $X$ be an object of a topos $\EE$ and for each $x\in X$ and $n\in\N$ define 
  $$t_{n, x}:=\tbigvee_{x\in A\in\mathcal{P}(X)^+} s_{n,A}\in \Delta_{\mathcal{P}(X)}^\EE.$$ \begin{enumerate}[label=\textup{(\arabic*)},leftmargin=2.0\parindent]
  \item \label{con.11} For each $n\in \N$, the localic map $q_n\colon Q_n:=\tbigcoprod^{\cat{Loc}}_{x\in X}\op(t_{n,x})\to \Delta_{\mathcal{P}(X)}^\EE$ is an open surjection;
  \item \label{con.12} The `multiple pullback' $q\colon Q\to  \Delta_{\mathcal{P}(X)}^\EE$ of $\{q_n\}_{n\in\N}$ is an open surjection;
  \item \label{con.13}The topos $\FF:=\text{Sh}_\EE(Q)$ contains a monomorphism $\N_\FF\rightarrowtail i^*(\mathcal{P}(X))$, where $i\colon \FF\hookrightarrow \EE$ is the inclusion.
  \end{enumerate}
\end{lemma}

\begin{proof}
\ref{con.11} follows at once from Lemma~\ref{lemma.open.coproducts} and the relation \ref{R1}.\\[2mm]
\ref{con.12} follows from \cite[Theorem~4.2]{Manuell2026}, which uses decidable equality of $\N$\footnote{We point out here that dedidable equality of the index set is essential contrarily to what is written in \cite[Theorem~3.5]{Plewe1997}.}, and the fact that open surjections regarded as objects  in the slice category $\cat{Loc}_\EE/\Delta_{\mathcal{P(X)}}^\EE$ correspond to positive overt locales internal to $\text{Sh}_\EE( \Delta_{\mathcal{P(X)}}^\EE )$.\\[2mm]
\ref{con.13} First,  for each $n\in\N$  the map  $h_n\colon \mathcal{P}(X)\to Q_n$ given by $$h_n(A):= ( \tbigvee^{\op(t_{n,x})} \{ 1 \mid x\in A\}   )_{x\in X}$$
is clearly a frame homomorphism. Composing it with the canonical frame homomorphism  $p_n\colon Q_n\to Q$, one obtains a frame homomorphism $g_n\colon\mathcal{P}(X)\to Q$. Consequently, the map $g\colon \mathcal{P}(X)\to \mathcal{P}(\N)\oplus Q$ given by
$$g(A):=\tbigvee_{n\in\N} \{n\}\oplus g_n(A)$$ is a frame homomorphism too. We claim that the frame homomorphism $$(g,\iota_2)\colon \mathcal{P}(X)\oplus Q\to \mathcal{P}(\N)\oplus Q$$ is surjective. First, we note that for each $n\in \N$, the relation $h_n(\{x\})\leq q_{n}^*(t_{n,x})$ holds by the definition of the maps $h_n$ and $q_n$. Then,  for $n\ne m$ and $x\in X$ one has
\begin{align*}
g_n(\{x\})\wedge g_m(\{x\})&=p_n(h_n(\{x\}))\wedge p_m(h_m(\{x\}))
\\& \leq p_nq_{n}^*(t_{n,x})\wedge p_mq_{m}^*(t_{m,x})
\\& = p_mq_{m}^*(t_{n,x})\wedge p_mq_{m}^*(t_{m,x})
\\& = p_mq_{m}^*(t_{n,x} \wedge t_{m,x})=0
\end{align*}
where the last equality holds because $t_{n,x} \wedge t_{m,x}=0$ as a consequence of  \ref{R3}.  Consequently, for each $n\in\N$ and $q\in Q$, one has
\begin{align*}
(g,\iota_2)(\tbigvee_{x\in X}\{x\}\oplus (g_{n}(\{x\})\wedge q)) 
&= \tbigvee_{x\in X} \tbigvee_{m\in \N} \{m\}\oplus g_m(\{x\})\wedge g_n(\{x\})\wedge q
\\&= \tbigvee_{x\in X} \{n\}\oplus g_n(\{x\})\wedge q=\{n\}\oplus q, 
\end{align*}
which shows that $(g,\iota_2)$ is surjective.  In other words, we have a diagram of localic maps
\[\begin{tikzcd}
	{\mathcal{P}(\N)\oplus Q} && {\mathcal{P}(X)\oplus Q} \\
	& Q
	\arrow["{(g,\iota_2)_*}", hook, from=1-1, to=1-3]
	\arrow[from=1-1, to=2-2]
	\arrow[from=1-3, to=2-2]
\end{tikzcd}\]
There is a well-known equivalence between the category $\mathbf{LH}_\EE/Q$ of local homeomorphisms over $Q$ and the topos $\text{Sh}_\EE(Q)$ (cf. \cite[Theorem~C.1.3.11]{elephant}). The unique locale map $!_Q\colon Q\to \Omega_\EE$ yields  by pullback (cf. \cite[Section~C.1.4]{elephant})  a functor $!_{Q}^*\colon \mathbf{LH}_\EE/\Omega_\EE \to \mathbf{LH}_\EE/Q$, which can be identified under the stated equivalence with  the left adjoint $i^* \colon \EE=\text{Sh}_\EE(\Omega_\EE)\to \text{Sh}_\EE(Q)$ of the inclusion $i\colon\text{Sh}_\EE(Q)\hookrightarrow \EE$. Moreover, under the equivalence $\EE\cong\mathbf{LH}_\EE/\Omega_\EE$, an object $X$ of $\EE$ corresponds to $\mathcal{P}(X)\to \Omega_\EE$. Consequently $i^*(X)$ (resp. $i^*(\N_\EE)$) are identified  under $\mathbf{LH}_\EE/Q\cong \text{Sh}_\EE(Q)$ with $!_{Q}^*(\mathcal{P}(X)\to \Omega_\EE)= \mathcal{P}(X)\oplus Q\to Q$ (resp. $\mathcal{P}(\N_\EE)\oplus Q\to Q$). Moreover, by virtue of \cite[Theorem~C.1.3.2\,(iii)]{elephant} the map $(g,\iota_2)_*$ is a monomorphism in $\mathbf{LH}_\EE/Q$, and so the equivalence $\mathbf{LH}_\EE/Q\cong \text{Sh}_\EE(Q)$ yields the desired monomorphism $i^*(\N_\EE)\rightarrowtail i^*(\mathcal{P}(X))$.
\end{proof}

At this point, we introduce the notion of geometric finiteness, due to Freyd (cf. \cite{Freyd2006}): 

\begin{definition}[\cite{Freyd2006}]
An object $X$ of a topos $\EE$ is \emph{geometrically finite} if for any geometric morphism $f\colon \mathcal{B}\to\mathcal{E}$ with $\mathcal{B}$ a Boolean topos, the object $f^*(X)$ is Kuratowski-finite.
\end{definition}

Geometric finiteness was studied by Henry in \cite{Henry2013}, where the following useful characterisation was proved:

\begin{lemma}[{\cite[Theorem~V.4, Corollary~V.7]{Henry2013}}]\label{charac.gf}
An object $X$ of a topos\footnote{We recall here that all toposes in this paper are assumed to have a natural numbers object.} $\EE$ is geometrically finite if and only if for any geometric morphism $f\colon \mathcal{B}\to\mathcal{E}$ with $\mathcal{B}$ a Boolean topos, if there is a monomorphism $\N_\mathcal{B}\rightarrowtail f^*(X)$ then $\mathcal{B}$ is degenerate.
\end{lemma}

We now prove the main result of this section:

\begin{theorem}\label{theorem.geometrically}
An object $X$ of a topos $\EE$ is geometrically finite if and only if the discrete locale $\mathcal{P}(X)$ is geometrically semi-irreducible in $\EE$.
\end{theorem}

\begin{proof}
Assume first that $X$ is geometrically finite. Consider the subtopos $i\colon  \text{Sh}_\EE(Q)\hookrightarrow \EE$ defined as in Lemma~\ref{lemma.QQN} where $m\colon \N_\FF\rightarrowtail i^*(\mathcal{P}(X))$ is a monomorphism. Consider the Booleanization $\mathcal{B}_Q$ of $Q$, which yields a Boolean topos $\mathcal{B}:=\text{Sh}_\EE(\mathcal{B}_Q)$ and a subtopos embedding $j\colon \text{Sh}_\EE(\mathcal{B}_Q)\hookrightarrow \text{Sh}_\EE(Q)$. Moreover, $j^*(m)\colon \N_\mathcal{B}\rightarrowtail (ij)^*(\mathcal{P}(X))$ is a monomorphism because $j^*$ preserves finite limits and, in particular, monomorphisms. Hence, by Lemma~\ref{charac.gf}, the topos $\text{Sh}_\EE(\mathcal{B}_Q)$ is degenerate---i.e., $\mathcal{B}_Q=\mathsf{O}$.  Since $\mathcal{B}_Q$ is dense in $Q$, it follows that $Q=\mathsf{O}$. By Lemma~\ref{lemma.QQN}\,\ref{con.12} there is a  surjection $Q\twoheadrightarrow \Delta_{\mathcal{P}(X)}^\EE$, so we conclude that $\Delta_{\mathcal{P}(X)}^\EE=\mathsf{O}$.

Suppose now that $\Delta_{\mathcal{P}(X)}^\EE=\mathsf{O}$. We show that $X$ is geometrically finite by checking the condition in Lemma~\ref{charac.gf}, so let $f\colon \mathcal{B}\to\mathcal{E}$ a geometric morphism with $\mathcal{B}$ a Boolean topos, and suppose there is a monomorphism $\N_\mathcal{B}\rightarrowtail f^*(X)$. By Lemma~\ref{first.lemma.geom}, the locale $f^\#( \Delta_{\mathcal{P}(X)}^\EE)$ has a point. But since $\Delta_{\mathcal{P}(X)}^\EE=\mathsf{O}$, it follows that $f^\#( \Delta_{\mathcal{P}(X)}^\EE)=\mathsf{O}$. The existence of a point means that there is a frame homomorphism $\mathsf{O}\to \Omega_\mathcal{B}$, which shows that $\Omega_\mathcal{B}$ is degenerate, and so $\mathcal{B}$ is degenerate too.
\end{proof}

\begin{corollary}\label{corl..fin}
An object $X$ of a topos $\EE$ is geometrically finite if and only if the discrete locale $\mathcal{P}(X)$ is $\FF$-compact.
\end{corollary}

\begin{remark}
The constructive Kuratowski--Mrówka theorem for locales, due to Vermeulen \cite{propermaps}, says that a locale $L$ is compact (in the sense that every open cover of $L$ has a Kuratowski-finite subcover) if and only if the product projection $\pi_2\colon L\oplus M\to M$ is closed for every locale $M$. Specialising to discrete locales, one has that a set $X$ is Kuratowski-finite if and only if the product projection $\pi_2\colon \mathcal{P}(X)\oplus M\to M$ is closed for every locale $M$. Hence Corollary~\ref{corl..fin} allows one to view geometric finiteness and Kuratowski-finiteness as dual properties in the sense that they correspond to categorical compactness with respect to the usual closure and the dual closure, respectively.
\end{remark}

\section{What happens classically}\label{section.whcla}

Throughout this section we shall freely use the law of excluded middle; in particular every locale is overt and an element is positive if and only if it is nonzero. We show that, classically,  $\FF$-compactness  is equivalent to semi-irreducibility in the sense of Dube \cite{irred} (cf. also the beginning of Section~\ref{sect4})  and moreover, we provide a simpler, direct proof of the characterisation.

We begin by showing that classically the locale of pairwise disjoint sequences has enough points.

\begin{proposition}\label{clasically.spatial}
Let $L$ be a locale and consider the discrete topological space $X_L:=(L^+,\mathcal{P}(L^+))$ on the set $L^+$ of nonzero elements of $L$.  
Then the frame $\Delta_{L}$ is isomorphic to the frame of opens of
$$Y_L:=\{ (a_n)_{n\in \N}\in X_{L}^{\N} \mid   n\ne m \implies a_n\wedge a_m=0\},$$
where $Y_L$ is considered as a subspace of  $X_{L}^{\N}$ (the latter endowed with the product topology). 
\end{proposition}

\begin{proof}
First, note that by the usual description of frame coproducts, the frame $\mathcal{P}(L^+)^\N$ (the $\N$-fold copower in the category of frames) can be presented by a generator $t_{n,a}$ for each $n\in\N$ and $a\in L^+$ subject to the relations $\tbigvee_{a\in L^+} t_{n,a}=1$ for every $n\in\N$ and $t_{n,a}\wedge t_{n,b}=0$ for every $n\in \N$ and $a\ne b$. Consequently, there is a frame homomorphism  $f\colon \mathcal{P}(L^+)^\N\to \Delta_L$  with $f(t_{n,a})=s_{n,a}$ because of \ref{R1} and \ref{R2}. Let $u:=\tbigvee_{a\wedge b>0,\,n\ne m} t_{n,a}\wedge t_{m,b}\in \mathcal{P}(L^+)^\N$. Because of \ref{R3}, one has $f(u)=0$ and so $f$ restricts to a frame homomorphism $f\colon \cl(u)\to \Delta_L$. Conversely, the assignment $s_{n,a}\mapsto t_{n,a}\vee u\in \cl(u)$ clearly sends \ref{R1}--\ref{R3} to identities; hence it defines a frame homomorphism $g\colon \Delta_L\to \cl(u)$. Now, $f$ and $g$ are inverse to each other so $\Delta_L\cong\cl(u)$.  It is well known that under classical logic countable products of discrete locales are spatial (cf. \cite[Theorem~10]{DowkerStrauss1977}), so the canonical frame homomorphism $\Omega(X_L)^\N \to \Omega (X_{L}^\N)$ is an isomorphism, and this clearly restricts to an isomorphism $\cl(u)\cong \Omega(Y_L)$.
\end{proof}

From the previous proposition and Corollary~\ref{MAINRESULT} one immediately obtains the following:

\begin{corollary}
A locale is $\FF$-compact if and only if it is semi-irreducible.
\end{corollary}

In the reminder of this section, we provide a simpler proof of the equivalence between $\FF$-compactness and semi-irreducibility under classical logic.

We begin by characterising the $\FF$-compact discrete locales.

\begin{lemma}{\label{prop.comp.disc}}
For a set $X$ of cardinality $\kappa$ and a locale $L$, denote  $\pi_2\colon\mathcal{P}(X)\oplus L\to L$ the product projection. Then the following are equivalent:
\begin{enumerate}[label=\textup{(\arabic*)},leftmargin=2.0\parindent]
\item\label{fit.cl.pr1} $\pi_2\in\FF$;
\item\label{fit.cl.pr2} $\kappa$-joins of fitted sublocales  of $L$ are fitted.
\end{enumerate}
\end{lemma}

\begin{proof}
Under the isomorphism $\mathcal{P}(X)\oplus L\cong \cat{Set}(X,L)$ (where the order is given pointwise), the projection is identified with the map $\pi \colon\cat{Set}(X,L)\to L$ given by $\pi(f)=\tbigwedge f[X]$. Then $\pi\in \FF$, iff for every family $\{ f_i\colon X\to L\}_{i\in I}$, the inclusion 
$$\tbigcap_{\pi[\tbigcap_i \op(f_i)]\subseteq \op(a)}\op(a)\subseteq \pi[\tbigcap_i \op(f_i)].$$
This is easily seen to be equivalent to the condition
$$\tbigcap_{\tbigcap_i \op(f_i)\subseteq \op(c_a)}\op(a)\subseteq \pi[\tbigcap_i \op(f_i)],$$
where $c_a$ denotes the constant function with value $a$. Thus the condition is equivalent to asking that for each $a,b\in L$ such that $\tbigcap_i \op(f_i)\subseteq \op(c_a)$ implies $b\in\op(a)$,  there is an $h\in \tbigcap_i\op(f_i)$ such that $b=\tbigwedge h[X]$. Since the Heyting operator in \cat{Set}$(X,L)$  is given pointwise, the latter is further equivalent to the condition that for each $a,b\in L$ such that for all $x\in X$,  $\tbigcap_i \op(f_i(x))\subseteq \op(a)$ implies $b\in\op(a)$, then there is a family $\{h_x\}_{x\in X}$ with $h_x\in \tbigcap_i\op(f_i(x))$ for each $x\in X$ such that $b=\tbigwedge_{x\in X} h_x$. The last part is equivalent to  $b\in \tbigvee_{x\in X} \tbigcap_{i\in I} \op(f_i(x))$ whereas the first part is equivalent to $b\in \left(\tbigvee_{x\in X} \tbigcap_{i\in I} \op(f_i(x))\right)^\circ$. Hence \ref{fit.cl.pr1} holds iff for every family $\{f_i\}_{i\in I}$, the sublocale $ \tbigvee_{x\in X} \tbigcap_{i\in I} \op(f_i(x))$ is fitted and this is clearly equivalent to \ref{fit.cl.pr2}. 
\end{proof}

Recall that a frame $L$ has \emph{property (H)} if for any $1\ne a\not\leq b$ in $L$ there is a $u\in L$ such that $u\not\leq a$ and $u^*\not\leq b$. It is well known (cf. \cite{separation}) that a $T_0$ space $X$ is Hausdorff if and only if $\Omega(X)$ has property (H). A sublocale of $L$ is called a \emph{one-point sublocale} if it is of the form $\bl(p)=\{1,p\}$ for a prime element $p$ of $L$ (see \cite[Section~III.10]{PP12}). In this context, we have the following result:

\begin{lemma}[{\cite[Proposition~3.10]{Arrieta2024}}]{\label{HimpF}}
In a Hausdorff locale, one-point sublocales are fitted.
\end{lemma}

\begin{remark}
In what follows, we shall need an example of a locale with a countable family of fitted sublocales whose join is not fitted. In order to construct one such example, choose a regular countable space $X$ with a dense subspace whose complement is also dense (e.g., the $X=\Q$ with the subspaces $A=\Z[\tfrac{1}{2}]\setminus \Z$ and $\Z[\tfrac{1}{3}]\setminus \Z$). The simple extension $E$ of $A$ in $X$ is strongly  Hausdorff ---in particular it satisfies property (H)--- but $\Omega(E)$ is not fit by \cite[Corollary~A.7]{Banaschewski1992}. Hence it contains a closed sublocale $\cl(U)$ which is not fitted. But  by Lemma~\ref{HimpF}, $\cl(U)=\tbigvee_{x\in X-U}\bl( X-\overline{\{x\}})$ is a countable join of fitted sublocales. 
\end{remark}

\begin{corollary}\label{class.px}
Let $X$ be a set. The discrete locale $\mathcal{P}(X)$ is  $\FF$-compact if and only if $X$ is finite.
\end{corollary}

\begin{proof}
For the `only if' part, assume that $X$ is $\FF$-compact and infinite. Select a locale $M$ with a countable join of fitted sublocales which is not fit (e.g., the one given by the previous remark). By Lemma~\ref{prop.comp.disc},   $\pi_2\colon \mathcal{P}(X)\oplus M\to M$ is not fitted, and so $\mathcal{P}(X)$ is not compact.  
The reverse implication follows by Lemma~\ref{prop.comp.disc} and the fact that fitted sublocales are closed under finite joins.
\end{proof}

\begin{lemma}\label{lemmaforbool}Every Boolean $\FF$-compact locale is finite.
\end{lemma}

\begin{proof}
Let $\{a_n\}_{n\in \N}$ be a pairwise disjoint family of nonzero elements of $B$. Set $a:=\tbigvee_{n\in\N}a_n$, and define a map $$f\colon \mathcal{P}(\N)\to \newdownarrow a$$ given by $f(A)=\tbigvee_{n\in A} a_n$ for every $A\subseteq \N$. Clearly, $f$ preserves arbitrary joins and the top element. Moreover, it also preserves binary meets by pairwise disjointness. Hence $f$ is a frame homomorphism. The fact that each $a_n$ is nonzero ensures that $f$ is injective. Note that, by Proposition~\ref{properties.categorical.compactness.fitted}\,\ref{properties.categorical.compactness.fitted1}, the locale  $\newdownarrow a$ is $\FF$-compact because it is (isomorphic to) an open, and hence fitted, sublocale of $B$.

 Hence there is a surjective localic map $f_*\colon\newdownarrow a \to \mathcal{P}(\N)$. Because the codomain is Boolean, $f_*$ is trivially open. Open surjections are pullback stable in \cat{Loc} (see \cite[Proposition~C.3.1.11]{elephant}), so by Proposition~\ref{properties.categorical.compactness.fitted}\,\ref{properties.categorical.compactness.fitted3}  one concludes that $\mathcal{P}(\N)$ is $\FF$-compact. This is a contradiction by Corollary~\ref{class.px}. Hence every pairwise disjoint family of nonzero elements of $B$ is finite,  and so $B$ is finite.
\end{proof}

\begin{remark}
Note that the above proof makes a heavy use of non-constructive principles and, in fact, openness of $f_*$  implies the law of excluded middle. Indeed, since Booleanness descends along open surjections, if $f_*$ were open, then $\mathcal{P}(\N)$ would be Boolean.  But, as is well known, discrete locales are never Boolean constructively, except for the trivial locale (indeed, take any proposition $p$ and form $A:=\{n\in\N \mid p\}$; then $0\in A$ or $0\in A^c$; the former case implies $p$ while the latter implies $\neg p$ because $A\cap A^c=\varnothing$).
\end{remark}

\begin{corollary}
Every  $\FF$-compact locale is semi-irreducible.
\end{corollary}

\begin{proof}
Let $L$ be a an $\FF$-compact locale. Since the Booleanization $\mathcal{B}_L$ is a fitted sublocale of $L$ (cf. \cite[Proposition~3.1]{Arrieta2025}), by Proposition~\ref{properties.categorical.compactness.fitted}\,\ref{properties.categorical.compactness.fitted1} it follows that $\mathcal{B}_L$  is $\FF$-compact, so by Lemma~\ref{lemmaforbool}, $\mathcal{B}_L$  is finite. The result then follows by \cite[Proposition~4.5]{irred}.
\end{proof}

To prove the converse, we first need a couple of lemmas:

\begin{lemma}{\label{finite.joins}}
A finite join of  $\FF$-compact sublocales of $L$ is $\FF$-compact.
\end{lemma}

\begin{proof}
Let $S_i\subseteq  L$, for $i=1,\dots,k$  be $\FF$-compact sublocales of $L$ with $\tbigvee_{i=1}^k S_i=L$. Let $M$ by any locale and $S$ a fitted sublocale of $L\oplus M$. Then for each $i=1,\dots, k$, one has that $S\cap (\pi_{1})_{-1}[S_i]$ is fitted in $(\pi_{1})_{-1}[S_i]$. Since $(\pi_{1})_{-1}[S_i]$ is (isomorphic to) $S_i\oplus M$, by $\FF$-compactness it follows that $\pi_2[  S\cap (\pi_{1})_{-1}[S_i] ]$ is fitted in $M$. Since finite joins of fitted sublocales are fitted, direct image commutes with joins and inverse image commutes with finite joins, it follows that $\pi_2[\tbigvee_{i=1}^k S\cap (\pi_{1})_{-1}[S_i]]=\pi_2[ S\cap (\pi_{1})_{-1}[\tbigvee_{i=1}^kS_i]]=\pi_2[S]$ is fitted in $M$.
\end{proof}

We recall here \cite{irred} that a frame is \emph{irreducible} if $a\wedge b=0$ implies either $a=0$ or $b=0$--- i.e., if and only if $\mathcal{B}_L=\{0,1\}$. In the following results, we use the explicit description of the binary coproduct of  frames as the family of \emph{saturated} downsets of their cartesian product. We refer the reader to \cite[Section~IV.4]{PP12} for details and adopt the notation used there.

\begin{lemma}{\label{lemma.irred}}
Let $L$ be a irreducible frame and $M$ any frame. Then, for any $b\in M$ and any  $D\in L\oplus M$ one has $$D\to (1\oplus b)= 1\oplus [ (\tbigvee_{(c,d)\in D,c\ne 0} d)\to b ].$$
\end{lemma}

\begin{proof}
One has $(x,y)\in D\to (1\oplus b)$ if and only if $(x\oplus y) \cap D \subseteq 1\oplus b$, that is, if and only if $\tbigvee_{(c,d)\in D}(x\wedge c\oplus y\wedge d)=\tbigvee_{(c,d)\in D,c\ne 0}(x\wedge c\oplus y\wedge d)\subseteq 1\oplus b $ --- i.e., if and only if for every $(c,d)\in D$ with $c\ne 0$, one has either $x\wedge c=0$ or $y\wedge d\leq b$. By irreducibility of $L$, the latter is equivalent to either $x=0$ or for every $(c,d)\in D$ with $c\ne 0$, one has $y\wedge d\leq b$. Clearly, this condition amounts to $(x,y)\in  1\oplus [ (\tbigvee_{(c,d)\in D,c\ne 0} d)\to b ]$.
\end{proof}

\begin{corollary}
Every semi-irreducible locale is $\FF$-compact.
\end{corollary}

\begin{proof}By \cite[Proposition~4.5]{irred}, a locale $L$ is semi-irreducible if and only if it is a join of a finite number of closed irreducible sublocales. Therefore, in view of Lemma~\ref{finite.joins}, it suffices to show the result for $L$ irreducible.
Let $\{D_i\}_{i\in I}\subseteq L\oplus M$  and $b\in \pi_2[\tbigcap_{i\in I}\op(D_i)]^\circ$. We claim that $1\oplus b\in \tbigcap_{i\in I}\op(D_i)$. Let $i\in I$. By Lemma~\ref{lemma.irred}, showing that  $1\oplus b\in \op (D_i)$ is equivalent to showing that $1\oplus [ (\tbigvee_{(c,d)\in D_i,c\ne 0} d)\to b ]= 1\oplus b$, that is, $(\tbigvee_{(c,d)\in D_i,c\ne 0} d)\to b=b$, or equivalently that $b\in \op(\tbigvee_{(c,d)\in D_i,c\ne 0} d)$. 
For the latter, observe that
\begin{align*}\tbigcap_{i\in I}\op(D_i)&\subseteq \op(D_i)= \tbigvee_{(c,d)\in D_i,c\ne 0} \op(c\oplus d)\subseteq \tbigvee_{(c,d)\in D_i,c\ne 0} \op(1\oplus d)= \op(1\oplus (\tbigvee_{(c,d)\in D_i,c\ne 0} d))\\&=(\pi_{2})_{-1}[\op(\tbigvee_{(c,d)\in D_i,c\ne 0} d)],
\end{align*}
which by adjunction is equivalent to $\pi_2[\tbigcap_{i\in I}\op(D_i)]\subseteq \op(\tbigvee_{(c,d)\in D_i,c\ne 0} d)$. Since $b\in \pi_2[\tbigcap_{i\in I}\op(D_i)]^\circ$, it follows that $b\in \op(\tbigvee_{(c,d)\in D_i,c\ne 0} d)$, as required. We therefore have $1\oplus b\in \tbigcap_{i\in I}\op(D_i)$ and so $b=\pi_2(1\oplus b)\in \tbigcap_{i\in I}\op(D_i)$.
\end{proof}

\noindent\textbf{Declaration of AI use.}
The author developed the  mathematical ideas and results presented in this paper. In preparing the manuscript, the author used ChatGPT (5.6 Sol and 6.1 Astra) to revise the exposition, polish the English, and help locate relevant literature. In particular, ChatGPT drew the author's attention to Freyd's notion of geometric finiteness~\cite{Freyd2006}.\\[2mm]
\noindent\textbf{Acknowledgement.}
The author acknowledges support from the Basque Government through grant IT1913-26, and a postdoctoral fellowship (grant POS-2022-1-0015 and POS-2025-2-0019).

\end{document}